\documentclass[numbers,webpdf,imanum]{ima-authoring-template}%

\graphicspath{{Fig/}}

\theoremstyle{thmstyletwo}%
\newtheorem{theorem}{Theorem}
\newtheorem{proposition}[theorem]{Proposition}%

\newtheorem{example}{Example}%
\newtheorem{remark}{Remark}%

\numberwithin{equation}{section}
\newtheorem{corollary}[theorem]{Corollary}%
\newtheorem{lemma}[theorem]{Lemma}%
\newtheorem{assumption}{Assumption}

\definecolor{darkgreen}{rgb}{0,0.5,0}
\usepackage{comment}
\usepackage{appendix}

\newcommand\Bregcheck{\Delta^{\!\!\check{\mathcal{R}}}}

\newcommand\ol[1]{\overline{#1}}

\newcommand\ul[1]{\underline{#1}}

\newcommand\res{\textrm{res}}
\newcommand\err{\textrm{err}}
\newcommand\sigmabg{\sigma_{bg}}
\newcommand\obsop{B}

\begin{document}

\DOI{DOI HERE}
\copyrightyear{2024}
\vol{00}
\pubyear{2024}
\access{Advance Access Publication Date: Day Month Year}
\appnotes{Paper}
\copyrightstatement{Published by Oxford University Press on behalf of the Institute of Mathematics and its Applications. All rights reserved.}
\firstpage{1}


\title[Convergence rates under a range invariance condition]{
Convergence rates under a range invariance condition with application to electrical impedance tomography}
\author{Barbara Kaltenbacher*
\address{\orgdiv{Department of Mathematics}, \orgname{University of Klagenfurt}, \orgaddress{\street{Universit\"atsstr. 65-67}, \postcode{9020}, \state{}, \country{Klagenfurt}}}}

\authormark{Barbara Kaltenbacher}

\corresp[*]{Corresponding author: \href{email:barbara.kaltenbacher@aau.at}{barbara.kaltenbacher@aau.at}}

\received{Date}{0}{Year}
\revised{Date}{0}{Year}
\accepted{Date}{0}{Year}


\abstract{
This paper is devoted to proving convergence rates of variational and iterative regularization methods under variational source conditions VSCs for inverse problems whose linearization satisfies a range invariance condition. In order to achieve this, often an appropriate relaxation of the problem needs to be found that is usually based on an augmentation of the set of unknowns and leads to a particularly structured reformulation of the inverse problem.
We analyze three approaches that make use of this structure, namely a variational and a Newton type scheme, whose convergence without rates has already been established in \cite{rangeinvar}; additionally we propose a split minimization approach that can be show to satisfy the same rates results.
\\
The range invariance condition has been verified for several coefficient identification problems for partial differential equations from boundary observations as relevant in a variety of tomographic imaging modalities.
Our motivation particularly comes from the by now classical inverse problem of electrical impedance tomography EIT and we study both the original formulation by a diffusion type equation and its reformulation as a Schr\"odinger equation. For both of them we find relaxations that can be proven to satisfy the range invariance condition.
Combining results on VSCs from \cite{Diss-Weidling} with the abstract framework for the three approaches mentioned above, we arrive at convergence rates results for the variational, split minimization and Newton type method in EIT. 
}


\maketitle

\section{Introduction}
Consider an inverse problem either in its
all-at-once formulation
\begin{equation}\label{aao}
\left.\begin{array}{ll}
A(q,u)=0 &\mbox{(model equation)}\\
\obsop u=y&\mbox{(observation equation)}
\end{array}\right\} 
\quad \Leftrightarrow: \ \mathbb{F}(q,u)=(0,y)^T
\end{equation}
or in its reduced formulation, with a parameter-to-state operator $S:\mathcal{D}(\mathbb{F})\to V$ that is implicitly defined by the first equation in 
\begin{equation}\label{red}
A(q,S(q))=0\mbox{ and }\obsop(S(q))=y \quad \Leftrightarrow: \ \mathbf{F}(q)=y.
\end{equation}
Here $A:\mathcal{D}(S)(\subseteq Q)\times V\to W^*$ is the (possibly nonlinear) model operator, $\obsop:V\to Z$ the (linear)  observation operator and $Q$, $V$, $W^*$, $Z$ are Banach spaces. 
Often $A$ represents the weak form of a partial differential equation (with boundary and -- if applicable -- initial conditions) and $W$ is the dual of a normed space.
More generally, consider
\begin{equation}\label{Fxy}
F(x)=y
\end{equation}
with an operator $F:\mathcal{D}(F)(\subseteq X)\to Y$, 
mapping between Banach spaces $X$ and $Y$, given noisy data $y^\delta$ with
\begin{equation}\label{delta}
\|y-y^\delta\|_Y\leq\delta
\end{equation}
and denoting by $x^\dagger$ an exact solution.
This comprises both all-at-once \eqref{aao} with $x=(q,u)$, $F=\mathbb{F}$ and reduced \eqref{red} with $x=q$, $F=\mathbf{F}$ formulations.
Since \eqref{Fxy} is typically ill-posed, regularization needs to be applied, cf., e.g., 
\cite{BakKok04, EnglHankeNeubauer:1996, Kirsch:2021, KalNeuSch08, SchusterKaltenbacherHofmannKazimierski:2012}. We will here consider variational and iterative approaches for this purpose.

\subsection*{Range invariance and other nonlinearity conditions}

In order to prove convergence for regularization methods, conditions on the nonlinearity of the forward operator $F$ are needed. One of the most commonly used of these is the {\em tangential cone} condition
\begin{equation}\label{tangcone}
\forall x,\tilde{x}\in U \,: \ \|F(x)-F(\tilde{x})- F'(x)(x-\tilde{x})\|_Y\leq c_{\textup{tc}}
\|F(x)-F(\tilde{x})\|_Y
\end{equation}
in a neighborhood $U$ of the exact solution. This condition (or closely related ones; cf. the introduction of \cite{rangeinvar} for a somewhat more detailed overview) allows to show convergence of 
iterative regularization methods \cite{HankeNeubauerScherzer:1995, DeEnSc98, IRGNMIvanov} as well as local convexity of the Tikhonov functional in variational regularization \cite{ChaventKunisch:1996}.
The control by the residual provided by \eqref{tangcone} allows to use the discrepancy principle as an a posteriori rule for the regularization parameter choice, which makes these methods - in particular the iterative ones - very convenient to implement.
However, the major drawback of \eqref{tangcone} is that it could so far only be verified in case of full observations or in a discretized (thus not ill-posed in the sense of instability) setting.
 
An alternative condition that avoids these limitations is {\em range invariance} of the linearized forward operator, cf., e.g., \cite{diss,NewtonKaczmarz,DeEnSc98,dpapertheo,KNSbook:2008,Broyden,SchHofNas23}, which we here (like in \cite{rangeinvar}) impose in a differential formulation:
\begin{equation}\label{rangeinvar_diff}
\exists x_0\in U\,, \ K\in L(\widetilde{X},Y)\, \forall x\in U \, \exists r(x)\in \widetilde{X}\,: \  F(x)-F(x_0)=Kr(x).
\end{equation}
This condition has been verified for a number of inverse problems with boundary (or actually arbitrary) observations, see, e.g., \cite{rangeinvar,UVA_nlIP,nonlinearity_imaging_both}.
However, as has been demonstrated there and we will also see in Examples~\ref{ex:Schroedinger} and ~\ref{ex:EIT} below, establishing \eqref{rangeinvar_diff}  sometimes requires extension of the original parameter space $\check{X}$ to some larger space $X$, that is, a relaxation of the original inverse problem, and this may lead to a loss of unique identifiability. Thus we add a penalty term that in the limit restricts reconstructions to the original parameter space, where they can (more likely) be shown to be unique. We do so by means of a penalty functional $\mathcal{P}:X\to[0,\infty]$ such that $\mathcal{P}(x^\dagger)=0$ for the exact solution $x^\dagger$, more concretely $\mathcal{P}(x^\dagger)=\|P(x)\|_X^p$ for some operator $P$ (typically a projection operator) and some $p\in[1,\infty)$. 

The inverse problem \eqref{Fxy} can therefore (locally in $U$) be rewritten as a system
\begin{equation}\label{FP}
\begin{aligned}
&K\hat{r}=y-F(x_0)&&\textup{ in }Y\\
&r(x)=\hat{r}&&\textup{ in }\widetilde{X}\\
&P(x)=0&&\textup{ in }X
\end{aligned}
\end{equation}
for the unknowns $(\hat{r},x)\in\widetilde{X}\times U$, $U\subseteq X$;
see also \cite[Remark 2.2]{DeEnSc98}.

\subsection*{Structure exploiting reconstruction approaches}

A natural way of applying {\em variational regularization} to \eqref{FP} is
\begin{equation}\label{var}
\begin{aligned}
&(\hat{r}_{\alpha,\beta}^\delta,x_{\alpha,\beta}^\delta)\in \mbox{argmin}_{(\hat{r},x)\in\widetilde{X}\times U} J_{\alpha,\beta}^\delta(\hat{r},x)\\
&\text{where } J_{\alpha,\beta}^\delta(\hat{r},x)
:=\underbrace{\|K\hat{r}+F(x_0)-y^\delta\|_Y^p+\alpha\mathcal{R}(\hat{r})}_{=:J_\alpha(\hat{r})}
+\underbrace{\beta\mathcal{Q}(r(x),\hat{r})^p + \|P(x)\|_X^p}_{=:J_\beta(\hat{r},x)} 
\end{aligned}
\end{equation}
with some $p\in[1,\infty)$, some 
regularization functional $\mathcal{R}:\widetilde{X}\to[0,\infty]$, and some distance functional 
$\mathcal{Q}:\widetilde{X}^2\to[0,\infty]$ such that $\mathcal{Q}(\hat{r}_1,\hat{r}_2)=0$ implies $\hat{r}_1=\hat{r}_2$.
Note that a minimizer of \eqref{var} can be computed more easily than for the original Tikhonov functional 
\begin{equation}\label{Tikh_check}
\check{J}_{\alpha}^\delta(\check{x})=\|Kr(\check{x})+F(x_0)-y^\delta\|_Y^p+\alpha\check{\mathcal{R}}(\check{x}),
\end{equation}
since if $\mathcal{R}$ is convex, so is $J_\alpha$ (in case of $p=2$ and $\mathcal{R}$ being defined by a squared norm, $J_\alpha$ is even quadratic) and for $\mathcal{Q}$ convex with respect to its first component, $J_\beta(\hat{r},\cdot)$ is locally uniformly convex. 

An efficient way to solve a structured minimization problem like  \eqref{var} is to alternatingly solve 
minimization problems for the two additive parts for a subset of the variables only and iterate this procedure, see, e.g. 
\cite{Auslender:1971,BertsekasTsitsiklis:1989,GrippoSciandrone:1999,OrtegaReinboldt:1970,LuoTseng:1992}. Since in our case $J_\alpha$ only depends on $\hat{r}$, this leads to the following (non-iterative) one-step {\em split minimization} procedure. 
\begin{equation}\label{sm}
\begin{aligned}
&\hat{r}_{sm}^\delta\in \mbox{argmin}_{\hat{r}\in\widetilde{X}} J_\alpha(\hat{r})
\text{ where }J_\alpha(\hat{r})=\|K\hat{r}+F(x_0)-y^\delta\|_Y^p+\alpha\mathcal{R}(\hat{r})\\
&x_{sm}^\delta\in \mbox{argmin}_{x\in U} J_\beta(\hat{r},x)
\text{ where }J_\beta(\hat{r},x)=\beta\mathcal{Q}(r(x),\hat{r}_{sm})^p + \|P(x)\|_X^p. 
\end{aligned}
\end{equation}

To define a regularized {\em (frozen) Newton} method, we locally linearize $r$ in \eqref{var}, assuming that it is G\^{a}teaux differentiable. 
\begin{equation}\label{Newton}
\begin{aligned}
&(\hat{r}_{n+1}^\delta,x_{n+1}^\delta)\in \mbox{argmin}_{(\hat{r},x)\in\widetilde{X}\times U} J_{n}^\delta(\hat{r},x)\\
&\mbox{where } J_{n}^\delta(\hat{r},x)
:=\|K\hat{r}+F(x_0)-y^\delta\|_Y^p+\alpha_n
\mathcal{R}(\hat{r})\\
&\hspace*{4cm}+\beta_n\mathcal{Q}(r(x_n^\delta)+r'(x_n^\delta)(x-x_n^\delta),\hat{r})^p + \|P(x)\|^p
\end{aligned}
\end{equation}
with some $p\in[1,\infty)$, $\alpha_n,\beta_n>0$ and functionals $\mathcal{Q}$, $\mathcal{R}$ as above.

In \cite{rangeinvar}, convergence without rates has been proven for \eqref{var}, \eqref{Newton} (but see Proposition~\ref{prop:conv} for the current, slightly different setting). 
In this paper we aim to establish convergence rates under additional regularity assumptions on the exact solution $x^\dagger$ that can be cast into so-called variational source conditions. Moreover, we will analyze convergence of the split minimization method \eqref{sm} under these conditions. 

\subsection*{Some illustrative examples}

\begin{example}\label{ex:Schroedinger}
Identify $c(x)$ (that is, $q=c$ in \eqref{aao} and \eqref{red}) in the elliptic PDE 
\begin{equation}\label{PDE_Schroedinger}
-\Delta u +c\, u=0\text{ in }\Omega  
\end{equation}
on a smooth bounded domain $\Omega\subseteq\mathbb{R}^d$, $d\in\{2,3\}$,
from the Neumann-to-Dirichlet N-t-D map $\Lambda\in L(H^{-1/2}(\partial\Omega),H^{1/2}(\partial\Omega))$; 
that is, $\check{\mathbf{F}}(c):=(\text{tr}^D_{\partial\Omega}u_{n})_{n\in\mathbb{N}}$ where $u_{n}$ solves \\
\begin{equation}\label{PDE_un}
-\Delta u_n +c\, u_n=0 \ \text{ in }\Omega, \quad 
\partial_\nu u_n=\varphi_n \ \text{ on }\partial\Omega, \quad \int_{\partial\Omega}u_{n}\, d\Gamma=0
\end{equation}
for a basis of boundary currents $\varphi_n\in L^2(\partial\Omega)$.
Here $\text{tr}^D_{\partial\Omega}:H^1(\Omega)\to H^{1/2}(\partial\Omega)$ denotes the Dirichlet trace operator; later on we will also use the Neumann trace 
$\text{tr}^N_{\partial\Omega}:H^s(\Omega)\to H^{s-3/2}(\partial\Omega)$, $s\in\mathbb{R}$. 

Note that the Hilbert-Schmidt norm of $\Lambda$ when considered as an operator from $L^2(\partial\Omega)$ into itself \cite[Section 6.1.2]{Diss-Weidling} is equivalent to the $l^2(\mathbb{N};L^2(\partial\Omega))$ norm of $\check{\mathbf{F}}(c)$, when ortho-normalizing the basis functions $\varphi_n$ in $L^2(\partial\Omega)$.

To establish \eqref{rangeinvar_diff}, we extend the parameter space 
\footnote{ the space $Q$ yet to be specified, cf. Example~\ref{ex:Schroedinger_details} below} 
by substituting $c\in Q$ with a sequence of potentials 
$\vec{c}=(c_j)_{j\in\mathbb{N}}\in\ell^\infty(Q)$ 
and on the other hand penalizing $\|P\vec{c}\|$, where  
\begin{equation}\label{P_Schroedinger}
P\vec{c}=\Bigl(c_j-\bigl(\sum_{k\in \mathbb{N}}w_k\bigr)^{-1} \, \sum_{k\in \mathbb{N}}w_k c_k\Bigr)_{j\in\mathbb{N}}.
\end{equation}
for some $\vec{w}\in \ell^1(\mathbb{R}^+)$ (e.g., $\vec{w}:=(j^{-2})_{j\in\mathbb{N}}$). 
The $\ell^\infty$ norm is used 
to allow for the true solution -- the constant sequence with value $c$ -- to be contained in this 
space. 
The projection $I-P$ on the subspace of constant sequences
is Lipschitz continuous, cf. \eqref{Lipschitz_P} below.
Correspondigly, we re-define the forward operator 
\begin{equation}\label{Fred_Schroedinger}
\mathbf{F}(\vec{c})_n:= (\text{tr}^D_{\partial\Omega}u_{n})_{n\in\mathbb{N}}
\text{ where $u_{n}=S(\vec{c})_n$ solves }
\left\{\begin{array}{rll}
-\Delta u_n +c_n\, u_n&=0 \quad &\text{ in }\Omega \\
\partial_\nu u_n&=\varphi_n &\text{ on }\partial\Omega \\ 
\int_{\partial\Omega}u_{n}\, d\Gamma&=0.
\end{array}\right.
\end{equation}
Indeed, it is readily checked that by setting 
\begin{equation}\label{r_Schroedinger}
r(\vec{c})_j:=(c_j-c_{0,j})\frac{u_j}{u_{0,j}}, \quad K:=F'(\vec{c}_0)
\end{equation}
where $u_j=S(\vec{c})_j$, $u_{0j}=S(\vec{c}_0)_j$, we formally obtain \eqref{rangeinvar_diff}. 
The problem of dividing by zero in \eqref{r_Schroedinger} in view of the last condition in \eqref{Fred_Schroedinger} can be avoided by using an all-at-once formulation, in which we define the forward operator by 
\begin{equation}\label{Faao_Schroedinger}
\mathbb{F}(\vec{c},\vec{u})_n:= \left(\begin{array}{c}
v\mapsto \int_\Omega(\nabla u_n\cdot\nabla v+c_n\,u_n\,v)\, dx -\int_{\partial\Omega}\varphi_n\, v\, d\Gamma\\
\int_{\partial\Omega}u_{n}\, d\Gamma\\
\text{tr}^D_{\partial\Omega}u_{n}
\end{array}\right)
\end{equation}
Also in this setting, the function defined by \eqref{r_Schroedinger} can be used to formally verify \eqref{rangeinvar_diff}, but $u_{0,n}$ does not need to be a PDE solution and in particular does not need to have vanishing average over the boundary any more. 
\end{example}

\begin{example}\label{ex:EIT}
Identify $\sigma(x)$ (that is, $q=\sigma$ in \eqref{aao} and \eqref{red}) in 
\begin{equation}\label{PDE_EIT}
-\nabla\cdot(\sigma\nabla u)=0\textup{ in }\Omega  
\end{equation}
on a smooth bounded domain $\Omega\subseteq\mathbb{R}^d$, $d\in\{2,3\}$,
from the N-t-D map $\Lambda\in L(H^{-1/2}(\partial\Omega),H^{1/2}(\partial\Omega))$;\\
that is, $\check{\mathbf{F}}(\sigma)=(\text{tr}^D_{\partial\Omega}u_{n})_{n\in\mathbb{N}}$ where $u_{n}$ solves \\
\[
-\nabla\cdot(\sigma\nabla u_n)=0 \ \text{ in }\Omega, \quad \partial_\nu u_n=\varphi_n, \ \text{ on }\partial\Omega, \quad \int_{\partial\Omega}u_{n}\, d\Gamma=0
\]
for a basis of boundary currents $\varphi_n\in L^2(\partial\Omega)$.
This is a prototype version of the well-known electrical impedance tomography EIT problem going back to Calder\'{o}n \cite{Calderon1980}, and extensively studied in the inverse problems literature, see, e.g., \cite[Chapter 12]{MuellerSiltanen2012} or \cite[Chapter 6]{Diss-Weidling} for further references.

Here, in order to establish \eqref{rangeinvar_diff}, we do the extension on the observations, by introducing artificial boundary voltages $\vec{z}=(z_n)_{n\in\mathbb{N}}$, thus extending the forward operator to 
\begin{equation}\label{Fred_EIT}
\mathbf{F}(\sigma,\vec{z})=(\text{tr}^D_{\partial\Omega}u_{n}+z_n)_{n\in\mathbb{N}}
\text{ where $u_{n}=S(\vec{c})_n$ solves }
\left\{\begin{array}{rll}
-\nabla\cdot(\sigma\nabla u_n)&=0 \quad &\text{ in }\Omega \\
\partial_\nu u_n&=\varphi_n &\text{ on }\partial\Omega \\ 
\int_{\partial\Omega}u_{n}\, d\Gamma&=0.
\end{array}\right.
\end{equation}
and penalizing $\|P(\sigma,\vec{z})\|$, where 
\begin{equation}\label{P_EIT}
P(\sigma,\vec{z})=\vec{z},
\end{equation}

This formally satisfies \eqref{rangeinvar_diff} with 
\begin{equation}\label{r_EIT}
\begin{aligned}
&r(\sigma,\vec{z})_j:=(\sigma-\sigma_0,\ z_j-z_{0,j}+\text{tr}^D_{\partial\Omega}v_j)
, \quad K(\ul{d\sigma},\ul{d\vec{z}}):=\check{\mathbf{F}}'(\sigma_0)\ul{d\sigma}+\ul{d\vec{z}},\\
&\text{ where $v_j$ solves }
\left\{\begin{array}{rll}
-\nabla\cdot(\sigma\nabla v_j)&=\nabla\cdot\bigl((\sigma-\sigma_0)\nabla (u_j-u_{0,j})\bigr)
 \quad &\text{ in }\Omega \\
\partial_\nu v_j&=0 &\text{ on }\partial\Omega \\ 
\int_{\partial\Omega}v_j\, d\Gamma&=0.
\end{array}\right.
\end{aligned}
\end{equation}

The relaxation we have done here is by far less extensive than the one from Example~\ref{ex:Schroedinger}: In view of the bijectivity of the operator mapping $\sigma$ to the sequence of boundary voltage data corresponding to currents $(\varphi_n)_{n\in\mathbb{N}}$ it can be viewed as equivalent to adding a second conductivity function as an artificial parameter. 

Note that this extension follows a general concept of achieving \eqref{rangeinvar_diff} by extension in data space, that we will briefly discuss in Section~\ref{sec:generalF}. 

A well-known relation to Example~\ref{ex:Schroedinger} that has often been used for uniqueness and stability proofs and also for establishing variational source conditions \eqref{vsc} in \cite{Diss-Weidling} is given by the identity 
\begin{equation}\label{cfromsigma}
c=\Phi(\sigma)=\frac{\Delta\sqrt{\sigma}}{\sqrt{\sigma}}.
\end{equation}
that allows to substitute $u$ solving \eqref{PDE_EIT} by $\tilde{u}:=\sqrt{\sigma}\,u$ solving \eqref{PDE_Schroedinger}.
Vice versa, since $\sigma$ is assumed to take the known (often also constant) background value $\sigmabg$ outside a ball of radius $\rho$ contained in $\Omega$
\begin{equation}\label{sigmabg}
\sigma(x)=\sigmabg(x)\quad x\in \mathbb{R}^d\setminus \mathcal{B_\rho}(0),\quad
\mathcal{B}_{\rho+\epsilon}(0)\subseteq \Omega
\end{equation}
for some $\rho$, $\epsilon>0$, we can also uniquely determine $\sigma$ from $c$ by means of the elliptic boundary value problem 
\begin{equation}\label{sigmafromc}
-\Delta\sqrt{\sigma}+c\sqrt{\sigma}=0\ \text{ in }\mathcal{B_\rho}(0), \quad \sqrt{\sigma}=\sqrt{\sigmabg}\text{ on }\partial\mathcal{B_\rho}(0).
\end{equation}
\end{example}

The precise function space setting for both examples will be provided in Section~\ref{sec:EIT} below, where due to the relations \eqref{cfromsigma}, \eqref{sigmafromc} we have also subsumed Example~\ref{ex:Schroedinger} under ``application to EIT''.

\subsection*{Notation and some convex analysis tools}
With the dual space $X^*$ of some normed space $X$, we denote the dual pairing by $\langle x^*,x\rangle_{X^*,X}$ for $x^*\in X^*$, $x\in X$.
\\
For some proper convex functional $f$ on $X$, the subdifferential  
$\partial f(x)=\{\xi\in X^*\,:\, \langle \xi,\tilde{x}-x\rangle_{X^*,X}\leq f(\tilde{x})-f(x)\text{ for all }\tilde{x}\in\mathrm{dom}(f)\}$ is known to be nonempty. 
\\
We will use the Bregman distance:
\begin{equation}\label{Bregman}
\Delta^f_{\xi_2}(x_1,x_2)=f(x_1)-f(x_2)-\langle x_2,x_1-x_2\rangle_{X^*,X}
\quad\text{for some selection }\xi_2\in\partial f(x_2).
\end{equation}
The convex conjugate of a real function $f:\mathbb{R}\to\mathbb{R}$ is defined by
\begin{equation}\label{convexconjugate}
f^*(s)=\sup_{t\in\mathbb{R}}(s\,t-f(t)).
\end{equation}
We denote the (smallest) Lipschitz constant of a Lipschitz continuous map $f$ by $L_f$. 

\subsection*{Variational source conditions}
In order to prove rates, we will make use of a variational source condition VSC in the original formulation (before  extension of the parameter space) 
\begin{equation}\label{checkFcheckxy}
\check{F}(\check{x})=y
\text{ with }\check{F}:\mathcal{D}(\check{F})(\subseteq\check{X})\to Y, \
\check{F}(\check{x})=F(\check{x})=Kr(\check{x})+F(x_0)
\end{equation} 
being defined on a subspace $\check{X}$ of $X$ such that $P\vert_{\check{X}}=0$ and endowed with a proper regularization functional $\check{\mathcal{R}}:\check{X}\to[0,\infty]$. 
The VSC reads as follows
\begin{equation}\label{vsc}
\begin{aligned}
\exists \check{\xi}\in\partial\check{\mathcal{R}}(x^\dagger)\subseteq\check{X}^*\, \forall \check{x}\in \mathcal{L}_{\ol{R}}^{\check{\mathcal{R}}}: 
\quad
-\langle \check{\xi}, \check{x}-x^\dagger\rangle
\leq b \Bregcheck_{\check{\xi}}(\check{x},x^\dagger)
+\psi(\|\check{F}(\check{x})-\check{F}(x^\dagger))\|^p)
\end{aligned}
\end{equation}
for some $b\in(0,1)$
(with $\mathcal{L}_{\ol{R}}^{\check{\mathcal{R}}}=\{\check{x}\in\check{X}\, : \, \check{\mathcal{R}}(\check{x})\leq \ol{R}\}$ for some $\ol{R}\geq \check{\mathcal{R}}(x^\dagger)$). 
We use VSCs in this original resticted setting because for some important inverse problems such as EIT, they are readily available from the literature, see for example \cite{Diss-Weidling}.
Here $\psi$ is supposed to be an index function, that is, 
\begin{equation}\label{psi_index}
\psi:[0,\infty)\to[0,\infty)\text{ is continuous and monotonically increasing with }\psi(0)=0,
\end{equation}
in particular a low rate one
\begin{equation}\label{lowrate}
\psi(t)\geq c_\psi t, \ t>0.
\end{equation}
Additionally, we extend $\psi$ by $-\infty$ for $t<0$ and assume that 
\begin{equation}\label{psi_concave}
\begin{aligned}
&-\psi\text{ is convex, }(-\psi)^* \text{ is an index function and }
t\mapsto t\, (-\psi)^*(-\tfrac{1}{t}) \text{ is strictly monotone}
\end{aligned}
\end{equation}
 
The most commonly used classes of such index functions are H\"older type and logarithm type, respectively:
\begin{equation}\label{psiHoelderlog}
(a)\ \psi(t)= t^\mu, \qquad (b)\ \psi(t)=(-\log(\min\{t,t_0\}))^{-\nu}
\end{equation}
for some $\mu\in(0,1]$, $\nu>0$, $t_0>0$, cf, e.g., \cite[Definition 2.2 and (2.3a), (2.3b)]{Diss-Weidling}. 
The convex conjugates of $-\psi$ are given by \cite[(2.25 a), (2.25 b)]{Diss-Weidling}
\[
(a)\ (-\psi)^*(-\tfrac{1}{t})= c_\mu t^{\mu/(1-\mu)}, \qquad (b)\ \psi(-\tfrac{1}{t})=(-\log(\min\{t,t_0\}))^{-\nu}(1+o(1))
\]
with $c_\mu=\mu^{\mu/(1-\mu)}-\mu^{1/(1-\mu)}$ and satisfy \eqref{psi_concave}.

\medskip 

Since we will partly use VSCs in a reduced setting for proving rates in an all-at-once setting, it is important to note that they are equivalent for the reduced and the all-at-once formulation in the following sense.
\begin{lemma}\label{lem:vsc_red-aao}
Condition \eqref{vsc} for $\check{F}=\mathbb{F}$ as in \eqref{aao} implies \eqref{vsc} for $\check{F}=\mathbf{F}$ as in \eqref{red} with the same regularization functional $\check{\mathcal{R}}$ (depending only on $q$).\\
If for any $q\in U$, the operator $A(q,\cdot)$ has a uniformly Lipschitz continous inverse and $\check{\mathcal{R}}$ only depends on $q$, then  
condition \eqref{vsc} for $\check{F}=\mathbf{F}$ as in \eqref{red} implies \eqref{vsc} for $\check{F}=\mathbb{F}$ as in \eqref{aao} with $u^\dagger=S(q^\dagger)$, the same regularization functional $\check{\mathcal{R}}$, and $\psi(t)$ replaced by $\psi(\bar{C}t)$ with $\bar{C}=2^{p-1}\max\{1,\|\obsop\|^p L_{A^{-1}}^p\}$.
\end{lemma}
{\em Proof.}\ 
The first assertion is trivial.
The second one follows from monotonicity of $\psi$ and the fact that for any fixed $q\in Q$ and any $u\in V$
\[
\begin{aligned}
&\|\mathbf{F}(q)-\mathbf{F}(q^\dagger)\|_Y^p
=\|\obsop (S(q)-S(q^\dagger))\|_Y^p\leq \Bigl(\|\obsop (u-u^\dagger)\|_Y + \|\obsop\|\, \|u-S(q)\|\Bigr)^p\\
&\leq \Bigl(\|\obsop (u-u^\dagger)\|_Y + \|\obsop\|\, L_{A^{-1}}\|A(q,u)-A(q,S(q))\|\Bigr)^p\\
&\leq 2^{p-1}\Bigl(\|\obsop (u-u^\dagger)\|_Y^p + \|\obsop\|^p\, L_{A^{-1}}^p\|A(q,u)\|^p\Bigr)
\leq \bar{C}\|\mathbb{F}(q,u)-\mathbb{F}(q^\dagger,u^\dagger)\|_Y^p,
\end{aligned}
\]
since $A(q,S(q))=0$.
\begin{flushright}
$\diamondsuit$
\end{flushright}

Since for fixed $q\in U$, the forward problem $A(q,\cdot)=0$ is typically well-posed, Lipschitz continuity of the inverse of $A(q,\cdot)$ is a natural condition to hold.\\
Note that also here we have the respective all-at-once and reduced formulations before extension of the parameter space in mind in order to make use of the results from \cite[Section 5]{Diss-Weidling} for Example~\ref{ex:Schroedinger}; see Example~\ref{ex:Schroedinger_details} below. 

\subsection*{Plan of the paper}
Most of this paper is devoted to proving convergence rates of the methods \eqref{var}, \eqref{sm}, \eqref{Newton} in Section~\ref{sec:convrates}. 
In Section ~\ref{sec:EIT} we apply these results to the EIT Examples~\ref{ex:Schroedinger}, \ref{ex:EIT}, thereby making use of results on variational source conditions from \cite{Diss-Weidling}.

\section{Convergence rates}\label{sec:convrates}
In this section we will prove convergence rates for the three regularization schemes \eqref{var}, \eqref{sm}, \eqref{Newton} in the relaxed formulation \eqref{FP} under a VSC \eqref{vsc} in the original formulation \eqref{checkFcheckxy}. 
To this end, throughout this section we will assume that 
\begin{equation}\label{xdagx0checkX}
\begin{aligned}
(I-P)(X)\subseteq \check{X}\,, \quad P(x^\dagger)=P(x_0)=0, \text{ thus }\quad x^\dagger, \ x_0 \in \check{X}. 
\end{aligned}
\end{equation} 
This is satisfied, e.g, if $X$ is a Hilbert space, $I-P=\text{Proj}_{\check{X}}$, (the orthogonal projection onto the subspace $\check{X}$). 

The key steps for deriving convergence rates for Tikhonov regularization \eqref{Tikh_check} in the original formulation \eqref{checkFcheckxy} are easily explained in case $\delta=0$. Firstly, minimality yields 
$\check{J}_\alpha^\delta(\check{x}_\alpha^\delta)\leq \check{J}_\alpha^\delta(x^\dagger)$ and thus by \eqref{delta}, \eqref{Bregman}, \eqref{vsc}
\[
\|\check{F}(\check{x}_\alpha^\delta)-\check{F}(x^\dagger)\|^p+(1-b)\alpha\Bregcheck(\check{x}_\alpha^\delta,x^\dagger)
\leq \psi(\|\check{F}(\check{x}_\alpha^\delta)-\check{F}(x^\dagger)\|^p),
\]
that is, an estimate in terms of the residual $\|\check{F}(\check{x}_\alpha^\delta)-\check{F}(x^\dagger)\|^p$ and the error $\Bregcheck(\check{x}_\alpha^\delta,x^\dagger)$. This can be resolved for the two single quantities by means of the following
Lemma, that can be found in several recent publications on convergence rates, e.g., as part of \cite[Theorem 2.20 (a)]{Diss-Weidling}, and whose proof we provide in the appendix, for the sake of completeness.
\begin{lemma}\label{lem:rate}
For any $\alpha$, $c_1$, $c_2$, $C_3$, $d$, $\err$, $\res \ \in[0,\infty)$ the relation
\begin{equation}\label{minLemma}
c_1\res  + c_2\alpha\, \err 
\leq d + \alpha\,\psi(C_3\res)
\end{equation}
implies the estimates
\begin{equation}\label{rate_Lemma}
\err\leq \frac{1}{c_2}\Bigl( \frac{d}{\alpha} + (-\psi)^*\Bigl(-\frac{\tilde{c}_1}{\alpha}\Bigr)\Bigr), \quad
\res \leq \frac{2}{c_1}\Bigl(d + \alpha(-\psi)^*\Bigl(-\frac{\tilde{c}_1}{2\alpha}\Bigr)\Bigr),
\end{equation}
where $(-\psi)^*(z)=\sup_{s\in\mathbb{R}}(z\,s-(-\psi)(s))$ is the convex conjugate of $-\psi$ and $\tilde{c}_1=\tfrac{c_1}{C_3}$.
\end{lemma}

Combined with a proper choice of $\alpha=\alpha(\delta)$ this also yields convergence rates with noisy data. In the relaxed setting \eqref{FP} and for methods \eqref{var}, \eqref{sm}, \eqref{Newton}, the estimates resulting from minimality become substantially more involved and have to take into account and control the interplay among the three components in \eqref{FP}.

\subsection{Variational regularization}\label{sec:var}
Existence of a minimizer of \eqref{var} can be established under certain coercivity, lower semicontinuity and closedness conditions on $J_{\alpha,\beta}^\delta$ and $\widetilde{X}\times U$, cf. \cite{rangeinvar}. However, we can avoid them (as well as exact computation of a minimizer) by defining $(\hat{r}_{\alpha,\beta,\eta}^\delta,x_{\alpha,\beta,\eta}^\delta)\in\widetilde{X}\times U$ to be any pair satisfying 
\begin{equation}\label{var_eta}
\forall (\hat{r},x)\in\widetilde{X}\times U\, : \ J_{\alpha,\beta}^\delta(\hat{r}_{\alpha,\beta,\eta}^\delta,x_{\alpha,\beta,\eta}^\delta) \leq J_{\alpha,\beta}^\delta(\hat{r},x)+\eta
\end{equation}
for some positive tolerance $\eta>0$ that is chosen in dependence of the noise level $\delta$, cf., e.g., \cite{NeubauerScherzer1990} as is of course also the regularization parameter $\alpha=\alpha(\delta)>0$.
For the latter purpose, as opposed to the case of a tangential cone condition \eqref{tangcone}, where the discrepancy principle is a natural a posteriori choice, we here remain with an a priori choice:
\begin{equation}\label{apriori}
\eta\leq C_\eta\delta^p, \quad 
\underline{\tau}\delta^p\leq \varphi^{-1}(\alpha):=\alpha(-\psi)^*\left(-\frac{1}{2\bar{C}\alpha}\right)\leq \overline{\tau}\delta^p
\end{equation}
cf. \eqref{psi_concave} with some sufficiently large $\bar{C}>0$ (cf. \eqref{Cbar}, \eqref{Cbar_sm} below)
and some fixed safety factors $0<\underline{\tau}<\ol{\tau}$. 
The parameter $\beta$ will be chosen as a fixed constant independent of $\delta$ in most of what follows (except for Proposition~\ref{prop:conv} where the choice is more general but includes constant $\beta$) and $C_\eta>0$ is a fixed constant as well.

We first of all aim to prove that under a variational source condition \eqref{vsc} the rates
\begin{equation}\label{rate}
\Bregcheck_{\check{\xi}} ((I-P)(r^{-1}(\hat{r}_{\alpha,\beta,\eta}^\delta)),x^\dagger)=O(\widetilde{\psi}(\delta^p)), \quad
\|K(\hat{r}_{\alpha,\beta,\eta}^\delta-r(x^\dagger))\|_Y=O(\delta).
\end{equation}
hold with $\widetilde{\psi}$ related to $\psi$ (ideally $\widetilde{\psi}=\psi$), cf. Theorem~\ref{thm:var} below.
Later on, from this and further estimates, we will derive convergence rates on $x_{\alpha,\beta,\eta}^\delta$, cf. Corollary~\ref{cor:var}.

To be able to make use of the VSC \eqref{vsc}, we will require that the regularization functional used in \eqref{var} is related to the functional appearing in \eqref{vsc} by 
\begin{equation}\label{RcheckR}
\mathcal{R}(\hat{r}):=\check{\mathcal{R}}((I-P)(r^{-1}(\hat{r})))
\end{equation}
and that the ties between $x$ and $\hat{r}$ imposed by the $\mathcal{Q}$ term in $J_{\alpha,\beta}^\delta$ are strong enough to control the out-of-$\check{X}$ parts  
\begin{equation}\label{PrQ}
\forall \hat{r}_1,\hat{r}_2\in X\, : \ \|P (r^{-1}(\hat{r}_1)) - P(r^{-1}(\hat{r}_2))\|_X\leq C_Q \mathcal{Q}(\hat{r}_1,\hat{r}_2).
\end{equation}
The latter holds, e.g., if $r:U\mapsto r(U)$ is bijective with Lipschitz continuous inverse, 
$\mathcal{Q}(\hat{r}_1,\hat{r}_2)=\|\hat{r}_1-\hat{r}_2\|_X$, 
and $P$ is Lipschitz continuous, for example, the projection operator on a Hilbert space.

\begin{theorem}\label{thm:var}
Under the assumption \eqref{xdagx0checkX}, let $F$ satisfy \eqref{rangeinvar_diff} with Lipschitz continuous $r$ and let the variational source condition \eqref{vsc} hold with $\psi$ satisfying \eqref{psi_index}, \eqref{lowrate}, \eqref{psi_concave}.
Moreover, assume that $\mathcal{R}$, $\mathcal{Q}$ are chosen according to \eqref{RcheckR}, \eqref{PrQ}.

Then for any family of data $(y^\delta)_{\delta\in(0,\bar{\delta}]}$ satisfying \eqref{delta} and regularization parameters $\alpha=\alpha(\delta)$, $\eta=\eta(\delta)$ satisfying \eqref{apriori}, any  corresponding family of regularized reconstructions according to \eqref{var_eta} satisfies the convergence rate \eqref{rate} with \eqref{psitilde},
as well as 
\begin{equation}\label{rates_further}
\|P(x_{\alpha,\beta,\eta}^\delta)\|=O(\delta),\quad 
\mathcal{Q}(r(x_{\alpha,\beta,\eta}^\delta), \hat{r}_{\alpha,\beta,\eta}^\delta)=O(\delta).
\end{equation} 
\end{theorem}

{\em Proof.}\ 
As usual in the analysis of variational regularization methods, we start with a minimality estimate. Comparing $J_{\alpha,\beta}^\delta(\hat{r}_{\alpha,\beta,\eta}^\delta,x_{\alpha,\beta,\eta}^\delta)$ with $J_{\alpha,\beta}^\delta(r(x^\dagger),x^\dagger)$, using $P(x^\dagger)=0$ and
\begin{equation}\label{triangle}
\|K\hat{r}_{\alpha,\beta,\eta}^\delta+F(x_0)-y^\delta\|_Y^p=
\|K(\hat{r}_{\alpha,\beta,\eta}^\delta-r(x^\dagger))+y-y^\delta\|_Y^p
\geq 2^{1-p}\|K(\hat{r}_{\alpha,\beta,\eta}^\delta-r(x^\dagger))\|_Y^p-\delta^p
\end{equation}
we obtain
\begin{equation}\label{min0_var}
\begin{aligned}
&2^{1-p} \|K(\hat{r}_{\alpha,\beta,\eta}^\delta-r(x^\dagger))\|_Y^p  
+ \alpha\, (\mathcal{R}(\hat{r}_{\alpha,\beta,\eta}^\delta) -\mathcal{R}(r(x^\dagger)))\\
&+ \beta \mathcal{Q}(r(x_{\alpha,\beta,\eta}^\delta), \hat{r}_{\alpha,\beta,\eta}^\delta)^p
+ \|P(x_{\alpha,\beta,\eta}^\delta)\|^p
\leq (2+C_\eta)\delta^p. 
\end{aligned}
\end{equation}
Now, using the VSC \eqref{vsc} with $\check{x}=(I-P)r^{-1}(\hat{r}_{\alpha,\beta,\eta}^\delta)$
together with \eqref{Bregman}, \eqref{RcheckR}
we conclude 
\[
\begin{aligned}
&\mathcal{R}(\hat{r}_{\alpha,\beta,\eta}^\delta) -\mathcal{R}(r(x^\dagger))
= \check{\mathcal{R}}((I-P)(r^{-1}(\hat{r}_{\alpha,\beta,\eta}^\delta)))-
\check{\mathcal{R}}(x^\dagger)\\
&=\Bregcheck_{\check{\xi}} ((I-P)(r^{-1}(\hat{r}_{\alpha,\beta,\eta}^\delta)),
x^\dagger)
+ \langle \check{\xi}, (I-P)(r^{-1}(\hat{r}_{\alpha,\beta,\eta}^\delta))-x^\dagger\rangle
\end{aligned}
\]
where 
\[
\begin{aligned}
&-\langle \check{\xi}, (I-P)(r^{-1}(\hat{r}_{\alpha,\beta,\eta}^\delta))-x^\dagger\rangle\\
&\leq 
b \Bregcheck_{\check{\xi}}((I-P)r^{-1}(\hat{r}_{\alpha,\beta,\eta}^\delta), x^\dagger)
+\psi\left(\|K \Bigl(r\bigl((I-P)r^{-1}(\hat{r}_{\alpha,\beta,\eta}^\delta)\bigr) -r(x^\dagger)\bigr)\|^p\right).
\end{aligned}
\]
Here, using \eqref{xdagx0checkX}, \eqref{PrQ}, and the monotonicity of $\psi$, we can further estimate 
\begin{equation}\label{estKterm_psi}
\begin{aligned}
&\psi\Bigl(\|K \bigl(r\bigl((I-P)r^{-1}(\hat{r}_{\alpha,\beta,\eta}^\delta)\bigr) -r(x^\dagger)\bigr)\|^p\Bigr)\\
&\leq 
\psi\Bigl(2^{p-1}\Bigl(\|K (\hat{r}_{\alpha,\beta,\eta}^\delta-r(x^\dagger))\|^p
+\|K\|^p L_r^p  \|P r^{-1}(\hat{r}_{\alpha,\beta,\eta}^\delta)\|^p\Bigr)\Bigr)\\
&\leq 
\psi\Bigl(2^{p-1}\|K (\hat{r}_{\alpha,\beta,\eta}^\delta-r(x^\dagger))\|^p
+2^{2p-2}\|K\|^p L_r^p  \Bigl(\|P(x_{\alpha,\beta,\eta}^\delta)\|^p + 
C_Q^p\mathcal{Q}(r(x_{\alpha,\beta,\eta}^\delta),\hat{r}_{\alpha,\beta,\eta}^\delta)^p\Bigr)\Bigr).
\end{aligned}
\end{equation}
Altogether we obtain
\begin{equation}\label{est_min}
\begin{aligned}
&2^{1-p} \|K(\hat{r}_{\alpha,\beta,\eta}^\delta-r(x^\dagger))\|_Y^p  
+ \beta \mathcal{Q}(r(x_{\alpha,\beta,\eta}^\delta), \hat{r}_{\alpha,\beta,\eta}^\delta)^p
+ \|P(x_{\alpha,\beta,\eta}^\delta)\|^p\\
&\hspace*{3cm}
+ \alpha\, (1-b) \Bregcheck_{\check{\xi}}((I-P)r^{-1}(\hat{r}_{\alpha,\beta,\eta}^\delta), x^\dagger)
\\
&\leq 
\psi\Bigl(2^{p-1}\|K (\hat{r}_{\alpha,\beta,\eta}^\delta-r(x^\dagger))\|^p
+2^{2p-2}\|K\|^p L_r^p  \Bigl( 
C_Q^p\mathcal{Q}(r(x_{\alpha,\beta,\eta}^\delta),\hat{r}_{\alpha,\beta,\eta}^\delta)^p
+\|P(x_{\alpha,\beta,\eta}^\delta)\|^p\Bigr)\Bigr)\\
&\qquad+(2+C_\eta)\delta^p. 
\end{aligned}
\end{equation}
Abbreviating 
\[
\begin{aligned}
\err&:=(1-b)\Bregcheck_{\check{\xi}} ((I-P)(r^{-1}(\hat{r}_{\alpha,\beta,\eta}^\delta)),x^\dagger)
,\\
\res &:=
2^{1-p}\|K \Bigl(\hat{r}_{\alpha,\beta,\eta}^\delta)-r(x^\dagger)\Bigr)\|^p
+ \beta \mathcal{Q}(r(x_{\alpha,\beta,\eta}^\delta),\hat{r}_{\alpha,\beta,\eta}^\delta)^p
+ \|P(x_{\alpha,\beta,\eta}^\delta)\|^p, 
\end{aligned}
\]
we obtain that
\begin{equation}\label{min}
\res  + \alpha\, \err 
\leq (2+C_\eta)\delta^p + \alpha\,\psi(\bar{C}\res)
\end{equation}
with
\begin{equation}\label{Cbar}
\bar{C}:= 2^{2p-2}\max\{1,\|K\|^p L_r^p,\tfrac{\|K\|^p L_r^p C_Q^p}{\beta}\}.
\end{equation}

Now an application of Lemma~\ref{lem:rate} to \eqref{min} with the regularization parameter choice \eqref{apriori} and \eqref{psi_concave} implies the rate 
\[
\begin{aligned}
&\Bregcheck_{\check{\xi}} ((I-P)(r^{-1}(\hat{r}_{\alpha,\beta,\eta}^\delta)),x^\dagger)
\leq \bar{\bar{C}}
(-\psi)^*\Bigl(-\frac{1}{\bar{C}\alpha}\Bigr)
\leq \bar{\bar{C}}
(-\psi)^*\Bigl(-\frac{1}{2\bar{C}\varphi(\overline{\tau}\delta^p)}\Bigr)\\
&\|K(\hat{r}_{\alpha,\beta,\eta}^\delta-r(x^\dagger))\|_Y^p
\leq 2^{p-1}\Bigl((2+C_\eta)\delta^p+
\alpha(-\psi)^*\Bigl(-\frac{1}{2\bar{C}\varphi(\overline{\tau}\delta^p)}\Bigr)\Bigr)
\leq 2^{p-1}(2+C_\eta+ \overline{\tau})\delta^p
\end{aligned}
\]
with $\bar{\bar{C}}=\frac{1}{1-b} \Bigl(\frac{2+C_\eta}{\underline{\tau}}+1\Bigr)$; 
that is, \eqref{rate} with 
\begin{equation}\label{psitilde}
\widetilde{\psi}(d)=(-\psi)^*\Bigl(-\frac{1}{2\bar{C}\varphi(\overline{\tau}d)}\Bigr), 
\qquad
\varphi^{-1}(t)=t\, (-\psi)^*(-\frac{1}{2\bar{C}t}),
\end{equation}
cf. \eqref{psi_concave}.
\begin{flushright}
$\diamondsuit$
\end{flushright}

\begin{remark}\label{rem:alphachoice}
Note that the ideal choice of $\alpha=\alpha(\delta)$ minimizing  
$\frac{d}{\alpha} + (-\psi)^*\Bigl(-\frac{1}{\bar{C}\alpha}\Bigr)$ with $d=(2+C_\eta)\delta^p$ yields the best possible rate cf. \cite[Chapter IV]{Diss-Weidling} and the references therein.
\[
\begin{aligned}
\err&\leq
\inf_\alpha\Bigl(\frac{d}{\alpha} + (-\psi)^*\Bigl(-\frac{1}{\bar{C}\alpha}\Bigr)\Bigr)
=\inf_\alpha
\Bigl( \left(-\bar{C}d\right)\left(-\frac{1}{\bar{C}\alpha}\right) + (-\psi)^*\Bigl(-\frac{1}{\bar{C}\alpha}\Bigr)\Bigr)\\
&=-\sup_s
\Bigl( \left(\bar{C}d\right)\, s - (-\psi)^*(s)\Bigr)
=-(-\psi)^{**}\left(\bar{C}d\right)
=\psi\left(\bar{C}d\right).
\end{aligned}
\]

On the other hand, in the logarithmic case \eqref{psiHoelderlog} (b), it is readily checked that the simple choice $\alpha(\delta)\sim\delta^{\tilde{p}}$ with $\tilde{p}\in(0,p]$ yields the same rate.
\end{remark}

A convergence rate of $x_{\alpha,\beta,\eta}^\delta$ to $x^\dagger$ can be obtained by combining the $\Bregcheck_\xi$, $\mathcal{Q}$ and $P$ terms.
To this end, consider the more particular setting that $\check{\mathcal{R}}$ is defined by a norm on an $s$-convex space $X$ 
and correspondingly $\mathcal{Q}$ is defined by the norm on $\widetilde{X}$ 
\begin{equation}\label{RQ_var}
\check{\mathcal{R}}(\check{x}):=\|\check{x}\|_{\check{X}}^s\, \quad
\Bregcheck_{\check{\xi}}(\check{x},x^\dagger)\geq c_s \|\check{x}-x^\dagger\|_{\check{X}}^s\, \quad
\mathcal{Q}(\hat{r}_1,\hat{r}_2):=\|\hat{r}_1-\hat{r}_2\|_{\widetilde{X}}
\end{equation}
for all $\hat{r}_1,\hat{r}_2\in \widetilde{X}$, $\check{x}\in \check{X}$.
Here $s\geq2$, cf. e.g. \cite[Example 2.47]{SchusterKaltenbacherHofmannKazimierski:2012}.

\begin{corollary}\label{cor:var}
Under the assumption \eqref{xdagx0checkX} with an $s$-convex space $\check{X}$, let $F$ satisfy \eqref{rangeinvar_diff} with Lipschitz continuous $r$, $r^{-1}$, $P$, and let the variational source condition \eqref{vsc} hold with $\psi$ satisfying \eqref{psi_index}, \eqref{lowrate}, \eqref{psi_concave}.
Moreover, assume that $\mathcal{R}$, $\mathcal{Q}$ are chosen according to \eqref{RcheckR}, \eqref{RQ_var}.

Then 
\[
\|x_{\alpha,\beta,\eta}^\delta-x^\dagger\|_X^s=O(\widetilde{\psi}(\delta^p)).
\]
\end{corollary}
{\em Proof.}\ 
By the triangle inequality we have 
\[
\begin{aligned}
&\|x_{\alpha,\beta,\eta}^\delta-x^\dagger\|_X\\
&\leq \|P(x_{\alpha,\beta,\eta}^\delta)\|_X
+\|(I-P)(x_{\alpha,\beta,\eta}^\delta)-(I-P)(r^{-1}(\hat{r}_{\alpha,\beta,\eta}^\delta))\|_X 
+ \|(I-P)(r^{-1}(\hat{r}_{\alpha,\beta,\eta}^\delta))-x^\dagger\|_X\\
&\leq O(\delta)+ O(\widetilde{\psi}(\delta^p))^{1/s},
\end{aligned}
\]
where we have estimated 
\[
\|(I-P)(x_{\alpha,\beta,\eta}^\delta)-(I-P)(r^{-1}(\hat{r}_{\alpha,\beta,\eta}^\delta))\|_X  \leq L_{I-P} L_{r^{-1}} \mathcal{Q}(r(x_{\alpha,\beta,\eta}^\delta),\hat{r}_{\alpha,\beta,\eta}^\delta).
\]
and used \eqref{rates_further}.
\begin{flushright}
$\diamondsuit$
\end{flushright}

\medskip

For the sake of completeness and to show that these results are compatible with convergence without rates in case a source condition is missing or unknown, we provide a result (whose proof is very similar to the one of \cite[Theorem 2.3]{rangeinvar}) in the setting under consideration here in the appendix, see Proposition~\ref{prop:conv} there.

%
\subsection{Split minimization}\label{sec:sm}
As in Section~\ref{sec:var}, in both subproblems of \eqref{sm} we minimize only up to a tolerance $\eta\leq C_\eta\delta^p$, 
\begin{equation}\label{sm_eta}
\begin{aligned}
&\forall \hat{r}\in\widetilde{X}\, : \ J_{\alpha}^\delta(\hat{r}_{sm}^\delta) \leq J_{\alpha}^\delta(\hat{r})+\eta\\
&\forall x\in U\, : \ J_{\beta}^\delta(\hat{r}_{sm}^\delta,x_{sm}^\delta) \leq J_{\beta}^\delta(\hat{r}_{sm}^\delta,x)+\eta
\end{aligned}
\end{equation}
and compare with the value at $\hat{r}=r(x^\dagger)$ and at $x=x^\dagger$ to obtain
\begin{equation}\label{est_min0_sm}
\begin{aligned}
&\|K\hat{r}_{sm}^\delta+F(x_0)-y^\delta\|_Y^p+\alpha\mathcal{R}(\hat{r}_{sm}^\delta)
\leq \alpha\mathcal{R}(r(x^\dagger))+(1+C_\eta)\delta^p\\
&\beta\mathcal{Q}(r(x_{sm}),\hat{r}_{sm})^p + \|P(x_{sm})\|_X^p 
\leq  \beta\mathcal{Q}(r(x^\dagger),\hat{r}_{sm})^p+C_\eta \delta^p
\end{aligned}
\end{equation}
Note that the additional right hand side term $\beta\mathcal{Q}(r(x^\dagger),\hat{r}_{sm})^p$ is due to suboptimality of split as compared to joint minimization and needs a dedicated estimate.
For this reason, instead of \eqref{RcheckR} we set
\begin{equation}\label{RcheckR_sm}
\begin{aligned} 
&\mathcal{R}(\hat{r}):=\check{\mathcal{R}}((I-P)(r^{-1}(\hat{r}))) + 
\mathcal{R}_P(P(r^{-1}(\hat{r}))) \\
&\text{ with }
\mathcal{R}_P(P(r^{-1}(\hat{r})))\geq\gamma\psi(C_P\|P(r^{-1}(\hat{r}))\|^p)
\end{aligned}\end{equation}
for sufficiently large constants $\gamma$, $C_P>0$
\begin{equation}\label{gammaCbar}
\gamma> C_\psi, \quad C_P\geq 2^{p-1}\|K\|^p L_r^p
\end{equation}
and make an assumption that allows to control $\beta\mathcal{Q}(r(x^\dagger),\hat{r}_{sm})^p$ by means of the regularzation term
\begin{equation}\label{Q_sm}
\forall\hat{r}\in \widetilde{X}: \quad\mathcal{Q}(r(x^\dagger),\hat{r})^p
\leq \tilde{C}_Q
\Bigl(\Bregcheck_{\check{\xi}}((I-P)r^{-1}(\hat{r}), x^\dagger)
+\psi\Bigl(C_P\|P(r^{-1}(\hat{r}))\|^p\Bigr).
\end{equation}
We also use an alternative estimate of the $\psi$ term as compared to \eqref{estKterm_psi} that is obtained by replacing the assumption \eqref{PrQ} by 
\begin{equation}\label{Cpsi}
\forall a,b \in\mathbb{R}^+: \quad\psi(a+b)\leq C_\psi(\psi(a)+\psi(b))
\end{equation}
for some $C_\psi>0$ (which is satisfied for all relevant cases of $\psi$, 
see the appendix).

\begin{theorem}\label{thm:var_sm}
Under the assumption \eqref{xdagx0checkX}, let $F$ satisfy \eqref{rangeinvar_diff} with Lipschitz continuous $r$ and let the variational source condition \eqref{vsc} hold with $\psi$ satisfying \eqref{psi_index}, \eqref{lowrate}, \eqref{psi_concave}, \eqref{Cpsi}.
Moreover, assume that $\mathcal{R}$, $\mathcal{Q}$ are chosen according to \eqref{RcheckR_sm}, \eqref{gammaCbar}, \eqref{Q_sm}.

Then for any family of data $(y^\delta)_{\delta\in(0,\bar{\delta}]}$ satisfying \eqref{delta} and regularization parameters $\alpha=\alpha(\delta)$, $\eta=\eta(\delta)$ satisfying \eqref{apriori},  a corresponding family of regularized reconstructions according to \eqref{var_eta} or \eqref{sm_eta} satisfies the convergence rate \eqref{rate} with \eqref{psitilde}
as well as 
\begin{equation}\label{rates_further_sm}
\begin{aligned}
&\|P(x_{sm})\|^p=O(\widetilde{\psi}(\delta^p)), \quad 
\mathcal{Q}(r(x_{sm}^\delta), \hat{r}_{sm}^\delta)^p=O(\widetilde{\psi}(\delta^p))\\
&\|P(r^{-1}(\hat{r}_{sm}^\delta))\|_X^p=O(\psi^{-1}(\widetilde{\psi}(\delta^p))).
\end{aligned}
\end{equation} 
\end{theorem}

{\em Proof.}\ 
Assumption \eqref{Cpsi} allows to bound the $\psi$ term arising from setting  $\check{x}:=(I-P)r^{-1}(\hat{r}_{sm}^\delta)$ in the VSC \eqref{vsc} as follows
\begin{equation}\label{estKterm_psi_sm}
\begin{aligned}
&\psi\Bigl(\|K (r\bigl((I-P)r^{-1}(\hat{r}_{sm}^\delta)\bigr) -r(x^\dagger))\|^p\Bigr)\\
&\leq C_\psi \Bigl(
\psi\Bigl(2^{p-1}\|K (\hat{r}_{sm}^\delta-r(x^\dagger))\|^p \Bigr)
+\psi\Bigl(2^{p-1}\|K\|^p L_r^p  \|P r^{-1}(\hat{r}_{sm}^\delta)\|^p\Bigr)\Bigr).
\end{aligned}
\end{equation}
From \eqref{est_min0_sm}, analogously to \eqref{est_min}, using \eqref{triangle}, \eqref{RcheckR_sm}, we obtain
\begin{equation}\label{est_min1_sm}
\begin{aligned}
&2^{1-p} \|K(\hat{r}_{sm}^\delta-r(x^\dagger))\|_Y^p  
+ \alpha\, (1-b) \Bregcheck_{\check{\xi}}((I-P)r^{-1}(\hat{r}_{sm}^\delta), x^\dagger)\\
&\hspace*{3cm}
+ \alpha\, \Bigl(\gamma\psi\Bigl(C_P\|P(r^{-1}(\hat{r}_{sm}))\|^p\Bigr) - C_\psi \psi\Bigl(2^{p-1}\|K\|^p L_r^p  \|P r^{-1}(\hat{r}_{sm}^\delta)\|^p\Bigr)
\\
&\hspace*{3cm}\leq 
\alpha C_\psi
\psi\Bigl(2^{p-1}\Bigl(\|K (\hat{r}_{sm}^\delta-r(x^\dagger))\|^p
\Bigr)+(2+C_\eta)\delta^p\\
&\beta\mathcal{Q}(r(x_{sm}),\hat{r}_{sm})^p + \|P(x_{sm})\|_X^p 
\leq  \beta\mathcal{Q}(r(x^\dagger),\hat{r}_{sm})^p +C_\eta \delta^p, 
\end{aligned}
\end{equation}
(where the second estimate is just a repetition of the second estimate in \eqref{est_min0_sm} for the convenience of the reader).
This suggests to choose $\gamma$, $C_P$ such that \eqref{gammaCbar} holds 
and motivates the assumption \eqref{Q_sm} on $\mathcal{Q}$.

With this, an additive combination of both estimates in \eqref{est_min1_sm} (after multiplying the second one with $\frac{\alpha}{\beta}\underline{c}$, where $\underline{c}=\frac{\min\{1-b,\gamma-C_\psi\}}{2\tilde{C}_Q}$)
yields
\begin{equation}\label{est_min2_sm}
\begin{aligned}
&2^{1-p} \|K(\hat{r}_{sm}^\delta-r(x^\dagger))\|_Y^p  
+ \alpha\, \tfrac{\min\{1-b,\gamma-C_\psi\}}{2} 
\Bigl(\Bregcheck_{\check{\xi}}((I-P)r^{-1}(\hat{r}_{sm}^\delta), x^\dagger)
+\psi\bigl(C_P\|P(r^{-1}(\hat{r}_{sm}))\|^p\bigr)\Bigr)\\
&\hspace*{3cm}+ \alpha \underline{c}\Bigl(\mathcal{Q}(r(x_{sm}),\hat{r}_{sm})^p + \tfrac{1}{\beta}\|P(x_{sm})\|_X^p\Bigr)\\
&\hspace*{3cm}\leq 
C_\psi
\psi\Bigl(2^{p-1}\Bigl(\|K (\hat{r}_{sm}^\delta-r(x^\dagger))\|^p
\Bigr)+(2+C_\eta+\tfrac{\alpha}{\beta}\underline{c} C_\eta)\delta^p. 
\end{aligned}
\end{equation}
Setting
\[
\begin{aligned}
\err&:=
\tfrac{\min\{1-b,\gamma-C_\psi\}}{2C_\psi} 
\Bigl(\Bregcheck_{\check{\xi}}((I-P)r^{-1}(\hat{r}_{sm}^\delta), x^\dagger)
+\psi\bigl(C_P\|P(r^{-1}(\hat{r}_{sm}))\|^p\bigr)\Bigr)\\
&\hspace*{3cm}+ \tfrac{\underline{c}}{C_\psi}\Bigl(\mathcal{Q}(r(x_{sm}),\hat{r}_{sm})^p + \tfrac{1}{\beta}\|P(x_{sm})\|_X^p\Bigr)
,\\
\res &:=
\tfrac{2^{1-p}}{C_\psi}\|K (\hat{r}_{sm}^\delta)-r(x^\dagger)\Bigr)\|^p
\end{aligned}
\]
we can rewrite this as
\begin{equation}\label{min_sm}
\res  + \alpha\, \err 
\leq \tfrac{2+C_\eta+\tfrac{\alpha}{\beta}\underline{c} C_\eta}{C_\psi}\delta^p + \alpha\,\psi(\bar{C}\res)
\end{equation}
with 
\begin{equation}\label{Cbar_sm}
\bar{C}:= 2^{2p-2}.
\end{equation}
Estimate \eqref{est_min0_sm} (with the right hand side in the second line replaced by $(1+C_\eta)\delta^p$) also holds true for the original simultaneous minimization version \eqref{var}, so that we can draw the same conclusions on it.

Applying Lemma~\ref{lem:rate} to \eqref{min_sm}, analogously to Theorem~\ref{thm:var} we obtain the assertion.
\begin{flushright}
$\diamondsuit$
\end{flushright}

\medskip

To obtain a convergence rate of $x_{sm}^\delta$ to $x^\dagger$ we again assume that $\check{\mathcal{R}}$ is defined by a norm on an $s$-convex space $X$, but in order to be compatible with \eqref{Q_sm}, make a slightly different assumption on $\mathcal{Q}$ as compared to \eqref{RQ_var}, namely: 
\begin{equation}\label{RQ_var_sm}
\begin{aligned}
&\check{\mathcal{R}}(x):=\|x\|_{X}^s\, \quad
\Bregcheck_{\check{\xi}}(x,x^\dagger)\geq c_s \|x-x^\dagger\|_X^s\, \\
&\overline{C}_Q\mathcal{Q}(\hat{r}_1,\hat{r}_2)^p
\geq
\max\Bigl\{
\|(I-P)r^{-1}(\hat{r}_2) - (I-P)r^{-1}(\hat{r}_1)\|_X^s,\\
&\hspace*{4cm}
c_\psi C_P\|P(r^{-1}(\hat{r}_2))-P(r^{-1}(\hat{r}_1))\|^p 
\Bigr\}, 
\end{aligned}
\end{equation}
for all $\hat{r}_1,\hat{r}_2\in \widetilde{X}$, $x\in X$, with $c_\psi$ as in \eqref{psi_concave}.

\begin{corollary}\label{cor:var_sm}
Under the assumption \eqref{xdagx0checkX} with an $s$-convex space $\check{X}$, let $F$ satisfy \eqref{rangeinvar_diff} with Lipschitz continuous $r$, $r^{-1}$, and let the variational source condition \eqref{vsc} hold with $\psi$ satisfying \eqref{psi_index}, \eqref{lowrate}, \eqref{psi_concave}, \eqref{Cpsi}.
Moreover, assume that $\mathcal{R}$, $\mathcal{Q}$ are chosen according to \eqref{RcheckR_sm}, \eqref{gammaCbar}, \eqref{Q_sm}, \eqref{RQ_var_sm}.

Then 
\[
\|x_{sm}^\delta-x^\dagger\|_X^s=O(\widetilde{\psi}(\delta^p)).
\]
\end{corollary}
{\em Proof.}\ 
By the triangle inequality we have 
\[
\begin{aligned}
&\|x_{sm}^\delta-x^\dagger\|_X\\
&\leq \|x_{sm}^\delta-r^{-1}(\hat{r}_{sm}^\delta)\|_X
+\|r^{-1}(\hat{r}_{sm}^\delta)-x^\dagger\|_X\\
&\leq \|(I-P)(x_{sm}^\delta)-(I-P)(r^{-1}(\hat{r}_{sm}^\delta))\|_X
+\|P(x_{sm}^\delta)-P(r^{-1}(\hat{r}_{sm}^\delta))\|_X\\
&+\|(I-P)(r^{-1}(\hat{r}_{sm}^\delta))-x^\dagger\|_X
+\|P(r^{-1}(\hat{r}_{sm}^\delta))\|_X\\
&\leq 
O(\widetilde{\psi}(\delta^p)^{1/s})
+O(\widetilde{\psi}(\delta^p)^{1/p})
+O(\psi^{-1}(\widetilde{\psi}(\delta^p)\bigr)^{1/p})
\end{aligned}
\]
where we have estimated (cf. \eqref{rates_further_sm})
\[
\begin{aligned}
&\|(I-P)(x_{sm}^\delta)-(I-P)(r^{-1}(\hat{r}_{sm}^\delta))\|_X
\leq 
\overline{C}_Q^{1/s}\mathcal{Q}(r(x_{sm}^\delta),\hat{r}_{sm}^\delta)^{p/s}
=O(\widetilde{\psi}(\delta^p)^{1/s})
\\
&\|P(x_{sm}^\delta)-P(r^{-1}(\hat{r}_{sm}^\delta))\|_X
\leq \left(\tfrac{\overline{C}_Q}{C_Pc_\psi}\right)^{1/p} \mathcal{Q}(r(x_{sm}^\delta),\hat{r}_{sm}^\delta)
=O(\widetilde{\psi}(\delta^p)^{1/p})
\\
&\|(I-P)(r^{-1}(\hat{r}_{sm}^\delta))-x^\dagger\|_X
\leq c_s^{-1/s}\Bregcheck_{\check{\xi}}((I-P)r^{-1}(\hat{r}_{sm}^\delta)),x^\dagger)^{1/s}
=O(\widetilde{\psi}(\delta^p)^{1/s})
\\
&\|P(r^{-1}(\hat{r}_{sm}^\delta))\|_X
=O(\psi^{-1}(\widetilde{\psi}(\delta^p))^{1/p}).
\end{aligned}
\]
\begin{flushright}
$\diamondsuit$
\end{flushright}

\begin{remark}\label{rem:mathcalR_P}
Since $\psi$ (or at least its decay rate) needs to be known to enable construction of a regularization functional satisfying \eqref{RcheckR_sm}, a convergence result without source condition like Proposition~\ref{prop:conv} would not make sense here.
\end{remark}

\subsection{Frozen Newton}\label{sec:Newton}
A minimizer to \eqref{Newton} exists if $U$ is weakly closed (e.g., $U$ is closed and convex) and $\mathcal{P}$ is weakly lower semicontinuous, but as in \eqref{var_eta}, \eqref{sm_eta} we can avoid imposing conditions for existence of a minimizer and computing it exactly by adding a tolerance $\eta_n>0$ that tends to zero as $n\to\infty$, thus defining 
the Newton step by
\begin{equation}\label{Newton_eta}
\forall (\hat{r},x)\in\widetilde{X}\times U\, : \ J_n^\delta(\hat{r}_{n+1}^\delta,x_{n+1}^\delta) \leq J_n^\delta(\hat{r},x)+\eta_n
\end{equation}

More precisely, we will choose the tolerance, the regularization parameters and the stopping index according to
\begin{equation}\label{alpha_eta_nstar}
\forall n\in\mathbb{N}:\quad \alpha_{n+1}<\alpha_n,
\quad
\eta_n\leq C_\eta \varphi^{-1}(\alpha_n), \qquad
n_*=\min\{n\in\mathbb{N}\, : \,  \alpha_n\leq\varphi(\tau\delta^p) \}.
\end{equation}
for fixed $C_\eta$, $\tau>0$,
where we additionally assume the growth condition
\begin{equation}\label{growthphipsi}
\forall n\in\mathbb{N}:\quad (-\psi)^*\Bigl(-\frac{1}{2\bar{C}\alpha_n}\Bigr)\leq
C_{(-\psi)^*} (-\psi)^*\Bigl(-\frac{1}{2\bar{C}\alpha_{n+1}}\Bigr)
\end{equation}
to hold for some $C_{(-\psi)^*}>1$.
Moreover, we set 
\begin{equation}\label{betaalpha}
\beta_n=c_{\beta\alpha}\alpha_n \text{ with }0<c_{\beta\alpha}\leq \tfrac{\min\{1-b,\gamma-C_\psi\}}{4\tilde{C}_Q},
\end{equation}
define $\mathcal{R}$ according to \eqref{RcheckR_sm} with \eqref{gammaCbar} and assume $\psi$ to satisfy \eqref{Cpsi}.

A further smoothness assumtion on $r$ will have to be formulated in terms of $\mathcal{Q}$ so that we can control the Taylor remainder estimates that will appear in the analysis. In order not cloud the procedure too much with technicalities we here restrict ourselves to the simple norm-based choice 
\begin{equation}\label{Qnorm}
\mathcal{Q}(\hat{r}_1,\hat{r}_2):=\|\hat{r}_1-\hat{r}_2\|_{\widetilde{X}}
\end{equation}
that automatically satisfies \eqref{Q_sm}, \eqref{RQ_var_sm}, and that we have already encountered earlier cf. \eqref{RQ_var}.
The assumption we make on $r$ 
\begin{equation}\label{tangcone_r}
\forall x_1,x_2\in U: \quad \|r(x_1)-r(x_2)-r'(x_2)(x_1-x_2)\|\leq c_{tcr} \|r(x_1)-r(x_2)\|
\end{equation}
with $c_{tcr}>0$ small enough
\begin{equation}\label{ctcr}
c_{tcr}<2^{2-2p} \, \max\left\{1,\, C_{(-\psi)^*}(1+2^{1-p})+1\right\}^{-1}
\end{equation}
resembles the tangential cone condition but now is naturally satisfied, given that $r$ is a local 
homeomorphism 
and we have constrained the minimization to a sufficiently small neighborhood $U$ of $x^\dagger$.

\begin{theorem}\label{thm:var_Newton}
Under the assumption \eqref{xdagx0checkX} with a $p$-convex space $\check{X}$, let $F$ satisfy \eqref{rangeinvar_diff} with Lipschitz continuous $r$ satisfying \eqref{tangcone_r}, \eqref{ctcr}, and let the variational source condition \eqref{vsc} hold with $\psi$ satisfying \eqref{psi_index}, \eqref{lowrate}, \eqref{psi_concave}, \eqref{Cpsi}.
Moreover, assume that $\mathcal{R}$, $\mathcal{Q}$ are chosen according to \eqref{RcheckR_sm}, \eqref{gammaCbar}, \eqref{Qnorm}.

Then for any family of data $(y^\delta)_{\delta\in(0,\bar{\delta}]}$ satisfying \eqref{delta} and regularization parameters $\alpha_n$, $\eta_n$, $n_*$ satisfying \eqref{alpha_eta_nstar}, \eqref{growthphipsi}, a corresponding family of regularized reconstructions according to \eqref{Newton_eta} satisfies the convergence rate 
\begin{equation}\label{rates_Newton}
\begin{aligned}
&\|\hat{r}_{n_*}-r(x^\dagger)\| = O(\widetilde{\psi}(\delta^p))^{1/p}, \qquad
\|r(x_{n_*}^\delta)-r(x^\dagger)\| = O(\widetilde{\psi}(\delta^p))^{1/p},\\
&\|K \Bigl(\hat{r}_{n_*}^\delta)-r(x^\dagger)\Bigr)\|_Y = O(\delta), \qquad
\|P(x_{n_*}^\delta)\| = O(\delta)
\end{aligned}
\end{equation} 
with $\tilde{\psi}$ as in \eqref{psitilde}.
\end{theorem}

{\em Proof.}\ 
As in the analysis of \eqref{var_eta}, and \eqref{sm_eta} we start with a minimality estimate. Comparing $J_n^\delta(\hat{r}_{n+1}^\delta,x_{n+1}^\delta)$ with $J_n^\delta(r(x^\dagger),x^\dagger)$, using $P(x^\dagger)=0$ and an analogon to \eqref{triangle}, we obtain the estimate 
(cf. \eqref{min0_var} with linearization in the first argument of $\mathcal{Q}$)
\[
\begin{aligned}
&2^{1-p} \|K(\hat{r}_{n+1}^\delta-r(x^\dagger))\|_Y^p  
+ \alpha_n\, \Bigl(\mathcal{R}(\hat{r}_{n+1}^\delta) -\mathcal{R}(r(x^\dagger))\Bigr)\\
&\qquad+ \beta_n \mathcal{Q}(r(x_{n+1}^\delta), \hat{r}_{n+1}^\delta)^p
+ \|P(x_{n+1}^\delta)\|^p\\
&\leq 2\delta^p+\eta_n+
\beta_n\mathcal{Q}(r(x_n^\delta)+r'(x_n^\delta)(x^\dagger-x_n^\delta),r(x^\dagger))^p
\\
&\qquad+\beta_n\max\Bigl\{0,\mathcal{Q}(r(x_{n+1}^\delta), \hat{r}_{n+1}^\delta)-\mathcal{Q}(r(x_n^\delta)+r'(x_n^\delta)(x_{n+1}^\delta-x_n^\delta),\hat{r}_{n+1}^\delta)\Bigr\}^p. 
\end{aligned}
\]
Making use of \eqref{RcheckR_sm}, \eqref{Qnorm}, \eqref{tangcone_r} as well as the variational source condition \eqref{vsc} together with \eqref{Cpsi} 
and 
\[
\|r(x_{n+1}^\delta)-\hat{r}_{n+1}^\delta\|^p\geq 2^{1-p}\|r(x_{n+1}^\delta)-r(x^\dagger)\|^p-\|\hat{r}_{n+1}^\delta-r(x^\dagger)\|^p
\]
\[
\|r(x_{n+1}^\delta)-r(x_{n}^\delta)\|^p\leq 
2^{p-1}\|r(x_{n+1}^\delta)-r(x^\dagger)\|^p+2^{p-1}\|r(x_{n}^\delta)-r(x^\dagger)\|^p
\]
we arrive at 
\begin{equation}\label{min0_fNew}
\begin{aligned}
&2^{1-p} \|K(\hat{r}_{n+1}^\delta-r(x^\dagger))\|_Y^p  \\
&\qquad
+ \alpha_n\, \tfrac{\min\{1-b,\gamma-C_\psi\}}{2} 
\Bigl(\Bregcheck_{\check{\xi}}((I-P)r^{-1}(\hat{r}_{n+1}^\delta), x^\dagger)
+\psi\bigl(C_P\|P(r^{-1}(\hat{r}_{n+1}))\|^p\bigr)\Bigr)\\
&\qquad
+ \beta_n (2^{1-p}-2^{p-1}c_{tcr})\|r(x_{n+1}^\delta)-r(x^\dagger)\|^p
+\|P(x_{n+1})\|^p\\
&
\leq \beta_n (1+2^{p-1})c_{tcr}\|r(x_n^\delta)-r(x^\dagger)\|^p
+\beta_n \|\hat{r}_{n+1}-r(x^\dagger)\|^p\\
&\qquad+\alpha_n C_\psi
\psi\Bigl(2^{p-1}\Bigl(\|K (\hat{r}_{n+1}^\delta-r(x^\dagger))\|^p
\Bigr)
+2\delta^p+\eta_n.
\end{aligned}
\end{equation}
Here under condition \eqref{Q_sm} we can control $\beta_n \|\hat{r}_{n+1}-r(x^\dagger)\|^p$ by means of 
the second term on the left hand side,
provided $\beta_n$ is chosen according to \eqref{betaalpha}.
Setting
\begin{equation}\label{err_res_Newton}
\begin{aligned}
\err_{n+1}&:= \tfrac{c_{\beta\alpha}}{C_\psi} \Bigl(
\|\hat{r}_{n+1}-r(x^\dagger)\|^p
+(2^{1-p}-2^{p-1}c_{tcr})\|r(x_{n+1}^\delta)-r(x^\dagger)\|^p\Bigr)
,\\
\res_{n+1}&:=
2^{1-p}\|K \Bigl(\hat{r}_{n+1}^\delta)-r(x^\dagger)\Bigr)\|^p
+ \|P(x_{n+1}^\delta)\|^p\\
d_n&:= \tfrac{c_{\beta\alpha}}{C_\psi} \alpha_n (1+2^{p-1})c_{tcr}\|r(x_n^\delta)-r(x^\dagger)\|^p
+\frac{2\delta^p+\eta_n}{C_\psi}
\end{aligned}
\end{equation}
and $\bar{C}$ as in \eqref{Cbar_sm}
we can rewrite this as
\[
\res_{n+1}  + \alpha_n\, \err_{n+1} 
\leq d_n + \alpha_n\,\psi(\bar{C}\res_{n+1}),
\]
thus by Lemma~\ref{lem:rate}
\begin{equation}\label{min_Newton}
\err_{n+1}\leq \frac{d_n}{\alpha_n} + (-\psi)^*\Bigl(-\frac{1}{\bar{C}\alpha_n}\Bigr), \quad
\res_{n+1} \leq 2\Bigl(d_n + \alpha_n(-\psi)^*\Bigl(-\frac{1}{2\bar{C}\alpha_n}\Bigr)\Bigr).
\end{equation}
As a consequence we get the recursion
\[
\begin{aligned}
\err_{n+1}\leq q\,\err_n +
\frac{2\delta^p+\eta_n}{C_\psi\,\alpha_n} + (-\psi)^*\Bigl(-\frac{1}{\bar{C}\alpha_n}\Bigr)
\leq q\,\err_n + \hat{C}\, (-\psi)^*\Bigl(-\frac{1}{2\bar{C}\alpha_n}\Bigr)
\end{aligned}
\]
with 
\[
q:=\frac{1+2^{p-1}}{2^{1-p}-2^{p-1}c_{tcr}}c_{tcr}\in(0,1/C_{(-\psi)^*})\subseteq(0,1),
\]
which holds under the smallness condition \eqref{ctcr} on $c_{tcr}$,
and $\hat{C}=\frac{2}{\tau C_\psi}+\frac{C_\eta}{C_\psi}+1$ for $n<n_*$ by \eqref{alpha_eta_nstar}.
Resolving this recursion and using \eqref{alpha_eta_nstar}, \eqref{growthphipsi}, we obtain
\[
\begin{aligned}
&\err_{n_*}\leq q^{n_*} \err_0 
+ \hat{C}\sum_{n=0}^{n_*-1} q^{n_*-n-1} \, (-\psi)^*\Bigl(-\frac{1}{2\bar{C}\alpha_n}\Bigr)\\
&\leq \Bigl((q C_{(-\psi)^*})^{n_*} \frac{\err_0}{(-\psi)^*\left(-\frac{1}{2\bar{C}\alpha_0}\right)}
+\frac{\hat{C}}{q(1-q C_{(-\psi)^*})} \Bigr)\, (-\psi)^*\Bigl(-\frac{1}{2\bar{C}\alpha_{n_*}}\Bigr)\\
&\leq \tilde{C}\, (-\psi)^*\Bigl(-\frac{1}{2\bar{C}\varphi(\tau\delta^p)}\Bigr)
\end{aligned}
\]
with $\tilde{C}= \frac{\err_0}{(-\psi)^*\left(-\frac{1}{2\bar{C}\alpha_0}\right)}
+\frac{C}{q(1-q C_{(-\psi)^*})}$.
 
This together with the second estimate in \eqref{min_Newton} and the definition \eqref{err_res_Newton} yields the convergence rates.
\begin{flushright}
$\diamondsuit$
\end{flushright}

\subsection{A canonical relaxation leading to range invariance}\label{sec:generalF}
A general way of extending an arbitrary inverse problem 
\[
\check{F}(\check{x})=y
\]
with $\check{F}:\mathcal{D}(\check{F})(\subseteq\check{X})\to Y$ to an operator equation \eqref{Fxy} satisfying the range invariance condition \eqref{rangeinvar_diff} is by defining
\[
X=\check{X}\times Y, \quad 
x=(\check{x},z), \quad 
F(x)= \check{F}(\check{x})+z, \quad
P(x)= z
\]
Indeed, it is readily checked that \eqref{rangeinvar_diff} formally holds with 
\[
r(x)=\left(\begin{array}{c}
\check{x}-\check{x}_0\\ 
z-z_0+\bigl(\check{F}(\check{x})-\check{F}(\check{x}_0)-\check{F}'(\check{x}_0)(\check{x}-\check{x}_0)\bigr)
\end{array}\right), \qquad 
K\underline{dx}=\check{F}'(\check{x}_0)\underline{d\check{x}}+\underline{dz}.
\]
A sufficient condition for the requirements of the convergence results in Section~\ref{sec:convrates}, namely local bijectivity of $r$ and Lipschitz continuity of $r$ and its inverse, is 
\[
\|r(x_1)-r(x_2) -(x_1-x_2)\|\leq c_r \|(x_1-x_2)\|, \qquad x_1,x_2\in U\times Y
\]
for some $c_r\in(0,1)$, which in the above setting amounts to a Taylor remainder estimate of the form 
\begin{equation}\label{Taylor}
\|\check{F}(\check{x}_1)-\check{F}(\check{x}_2)-\check{F}'(\check{x}_0)(\check{x}_1-\check{x}_2)\|_Y
\leq c_r\|\check{x}_1-\check{x}_2\| \qquad \check{x}_1,\check{x}_2\in U.
\end{equation}
Likewise, the additional condition \eqref{tangcone_r} required for proving convergence rates for \eqref{Newton_eta} reads as 
\begin{equation}\label{tcr_Fcheck}
\|\check{F}(\check{x}_1)-\check{F}(\check{x}_2)-\check{F}'(\check{x}_2)(\check{x}_1-\check{x}_2)\|_Y
\leq c_{tcr}\|\check{x}_1-\check{x}_2\| \qquad \check{x}_1,\check{x}_2\in U.
\end{equation}

The methods analyzed above in this setting read as follows.
\begin{itemize}
\item
variational regularization \eqref{var}: 
\begin{equation}\label{var_genF}
\begin{aligned}
&(\hat{r}_{\check{x},\alpha,\beta}^\delta,\hat{r}_{z,\alpha,\beta}^\delta,\check{x}_{\alpha,\beta}^\delta,z_{\alpha,\beta}^\delta)\in \mbox{argmin}_{(\hat{r}_{\check{x}},\hat{r}_z,\check{x},z)\in\check{X}\times Y\times U\times Y} J_{\alpha,\beta}^\delta(\hat{r}_{\check{x}},\hat{r}_z,\check{x},z)\\
&\text{where } J_{\alpha,\beta}^\delta(\hat{r}_{\check{x}},\hat{r}_z,\check{x},z)=
J_{\alpha}(\hat{r}_{\check{x}},\hat{r}_z)
+J_{\beta}(\hat{r}_{\check{x}},\hat{r}_z,\check{x},z)\\
&J_{\alpha}(\hat{r}_{\check{x}},\hat{r}_z)=
\|\check{F}'(\check{x}_0)\hat{r}_{\check{x}}+\check{F}(\check{x}_0)+\hat{r}_z+z_0-y^\delta\|_Y^p
+\alpha\check{\mathcal{R}}(\hat{r}_{\check{x}})\\
&J_{\beta}(\hat{r}_{\check{x}},\hat{r}_z,\check{x},z)=
\beta\Bigl(\mathcal{Q}_{\check{x}}(\check{x}-\check{x}_0,\hat{r}_{\check{x}})^p \\
&\hspace*{3cm}+\mathcal{Q}_{z}(z-z_0+\bigl(\check{F}(\check{x})-\check{F}(\check{x}_0)-\check{F}'(\check{x}_0)(\check{x}-\check{x}_0)\bigr),\hat{r}_z)^p\Bigr) \\
&\hspace*{3cm}+ \|z\|_X^p 
\end{aligned}
\end{equation}
\item
split minimization \eqref{sm}: 
\begin{equation}\label{sm_genF}
\begin{aligned}
&(\hat{r}_{\check{x},sm}^\delta,\hat{r}_{z,sm}^\delta)\in \mbox{argmin}_{(\hat{r}_{\check{x}},\hat{r}_z)\in\widetilde{X}} J_\alpha(\hat{r}_{\check{x}},\hat{r}_z)
\\
&(\check{x}_{sm}^\delta,z_{sm}^\delta)\in \mbox{argmin}_{(\check{x},z)\in U\times Y} J_\beta(\hat{r}_{\check{x},sm}^\delta,\hat{r}_{z,sm}^\delta,\check{x},z)
\end{aligned}
\end{equation}
\item
frozen Newton \eqref{Newton}:
\begin{equation}\label{Newton_genF}
\begin{aligned}
&(\hat{r}_{\check{x},n+1}^\delta,\hat{r}_{z,n+1}^\delta,\check{x}_{n+1}^\delta,z_{n+1}^\delta)\in \mbox{argmin}_{(\hat{r}_{\check{x}},\hat{r}_z,\check{x},z)\in\check{X}\times Y\times U\times Y} J_{n}^\delta(\hat{r}_{\check{x}},\hat{r}_z,\check{x},z)\\
&\text{where } J_{n}^\delta(\hat{r}_{\check{x}},\hat{r}_z,\check{x},z)
:=J_{\alpha_n}(\hat{r}_{\check{x}},\hat{r}_z)+J_{\beta_n,n}(\hat{r}_{\check{x}},\hat{r}_z,\check{x},z)\\
&J_{\beta,n}(\hat{r}_{\check{x}},\hat{r}_z,\check{x},z)=
\beta\Bigl(\mathcal{Q}_{\check{x}}(\check{x}-\check{x}_0,\hat{r}_{\check{x}})^p \\
&\hspace*{3cm}+\mathcal{Q}_{z}(z-z_0+
\bigl((\check{F}'(\check{x}_n)-\check{F}(\check{x}_0))(\check{x}-\check{x}_n)+y_n\bigr)
,\hat{r}_z)^p\Bigr) \\
&\hspace*{3cm}+ \|z\|_X^p \\
&y_n=\check{F}(\check{x}_n)-\check{F}(\check{x}_0)-\check{F}'(\check{x}_0)(\check{x}_n-\check{x}_0)
\end{aligned}
\end{equation}
\end{itemize}

\section{Application to EIT}\label{sec:EIT}

\begin{example}[Example~\ref{ex:Schroedinger} revisited]\label{ex:Schroedinger_details}
We adopt a convenient function space setting for this problem from \cite[Section 6.2]{Diss-Weidling} and show that 
formal range invariance \eqref{rangeinvar_diff} with \eqref{r_Schroedinger} can be made rigorous in this framework. 
This is done 
in the all-at-once formulation, since the latter allows to avoid division by zero in \eqref{r_Schroedinger} by a proper choice of $u_0$.
A logarithmic variational source condition for the reduced formulation is taken from \cite{Diss-Weidling}, cf. Theorem~\ref{thm:vsc_Schroedinger} below, and extended to the all-at-once setting by means of Lemma~\ref{lem:vsc_red-aao}.

The inverse problem can be written in reduced form as 
\begin{equation}\label{red_Schroedinger}
\check{\mathbf{F}}(c)=y
\end{equation}
or equivalently in all-at-once form as
\begin{equation}\label{aao_Schroedinger}
\check{\mathbb{F}}(c,u):=
\left(\begin{array}{c}\check{A}(c,u)\\\obsop u\end{array}\right)=\left(\begin{array}{c}0\\y\end{array}\right)
\end{equation}
(note that we are using the $\check{\cdot}$ notation for the formulation before extension)
or equivalently, after extension of the parameter space and penalization, as
\begin{equation}\label{aao_Schroedinger_extended}
\begin{aligned}
&\mathbb{F}(\vec{c},u):=
\left(\begin{array}{c}A(\vec{c},u)\\\obsop u\end{array}\right)=\left(\begin{array}{c}0\\y\end{array}\right)\\
&P(\vec{c},u)=0
\end{aligned}
\end{equation}
cf. \eqref{aao} with 
\begin{equation}\label{DFcheck}
\begin{aligned}
&\check{\mathbf{F}}:\mathcal{D}(\check{\mathbf{F}})(\subseteq \check{Q})\to Z, \quad
\check{A}:\mathcal{D}(\check{A})(\subseteq\check{X})\to (0,\vec{\varphi},0)+W^*, \\ 
&A:\mathcal{D}(A)(\subseteq X)\to (0,\vec{\varphi},0)+W^*, \quad
\obsop:V\to Z\\
&\mathcal{D}(\check{\mathbf{F}})= \Phi(\mathcal{D}), \quad 
\mathcal{D}(\check{A})= \mathcal{D}(\check{\mathbf{F}})\times V, \quad
\mathcal{D}(A)= \ell^\infty(\mathcal{D}(\check{\mathbf{F}}))\times V\\
&\text{ where }
\mathcal{D}=\{\sigma\in H^{2+s_1}(\Omega)\, : \, 0<\underline{\sigma}\leq \sigma(x)\leq \overline{\sigma}<\infty \ x\in\Omega \text{ and \eqref{sigmabg} holds}\},
\\
&\check{Q}=H^{s_1}(\Omega), \quad 
Q=\ell^\infty(\check{Q}), \quad 
V=\ell^2(H^2(\Omega)), \quad
W^*=\ell^2(L^2(\Omega))\times\ell^2(H^{1/2}(\partial\Omega))\times\ell^2(\mathbb{R}), \quad 
\\
&\check{X}=\check{Q}\times V, \quad 
X=Q\times V, \quad 
Z=\ell^2(L^2(\partial\Omega)),\quad 
Y:=W^*\times Z.
\end{aligned}
\end{equation}
Note that taking $Z=\ell^2(L^2(\Omega))$ together with the $L^2(\partial\Omega)$ normalization of the currents $\varphi_n$ corresponds to taking the Hilbert-Schmidt norm of the N-t-D operator.
Also note that regardless of the use of $\ell^\infty$ in the definition of $X$, the original space $\check{X}$ is still $s$-convex with $s=2$.  
The operators are defined by 
\begin{equation}\label{A_Schroedinger}
\begin{aligned}
&A(\vec{c},u)=
\left(\begin{array}{c}
-\Delta u_n +c_n u_n\\
\text{tr}^N_{\partial\Omega} u_n-\varphi_n\\
\int_{\partial\Omega}u_{n}\, d\Gamma
\end{array}\right)_{n\in\mathbb{N}}\!\!\!\!\!\!\!\!, \quad
\check{A}(c,u)= A(\vec{c},u) \text{ with $c_n\equiv c$}, \\
&\obsop u=\left(\text{tr}^D_{\partial\Omega}u_n\right)_{n\in\mathbb{N}}\quad
P(\vec{c},u)=\Bigl(c_n-\bigl(\sum_{k\in \mathbb{N}}w_k\bigr)^{-1} \, \sum_{k\in \mathbb{N}}w_k c_k\Bigr)_{n\in\mathbb{N}},\\
&\check{\mathbf{F}}(c) = \mathbf{F}(c):= (\text{tr}^D_{\partial\Omega}u_{n})_{n\in\mathbb{N}}
\text{ where $u_{n}=S(c)_n$ and }\check{A}(c,S(c))=0
\\
&\mathbf{F}(\vec{c}):= (\text{tr}^D_{\partial\Omega}u_{n})_{n\in\mathbb{N}}
\text{ where $u_{n}=\vec{S}(\vec{c})_n$ and }A(\vec{c},\vec{S}(\vec{c}))=0
\end{aligned}
\end{equation}
with the Neumann trace operator $\text{tr}^N_{\partial\Omega}$, 
Note that slightly deviating from \cite{Diss-Weidling}, in order to avoid division by zero in \eqref{r_Schroedinger}, we use all of $H^2(\Omega)$ rather than incorporating the vanishing boundary average constraint into the state space.
Moreover, we increase regularity as compared to the $H^1(\Omega)$ setting from \cite{Diss-Weidling} as this is possible with $s_1\geq0$ (note that even $s_1\in(1/2,3/2)$ is assumed in Theorem~\ref{thm:vsc_Schroedinger}) and will be needed in the rigorous verification of \eqref{rangeinvar_diff}.

Indeed, the mapping property $A(\vec{c},u):\mathcal{D}(A)\to W^*$ holds true due to 
boundedness of $\text{tr}^N_{\partial\Omega}:H^2(\Omega)\to H^{1/2}(\partial\Omega)$, as well as continuity of the multiplication operator 
$M:H^{s_1}(\Omega)\times H^2(\Omega)\mapsto L^2(\Omega)$, $(c,u)\mapsto c\,u$ 
for $s_1\geq0$, 
that 
allows to estimate
\[
\begin{aligned}
&\|A(\vec{c},u)+(0,\vec{\varphi},0)\|_{W^*}\\
&\leq \Bigl(\sum_{n\in\mathbb{N}} 
\bigl(1+C_{H^2,L^\infty}^\Omega \|c_n\|_{L^2(\Omega)}
+\|\text{tr}^N_{\partial\Omega}\|_{H^2(\Omega),H^{1/2}(\partial\Omega)}
+\|\text{tr}^D_{\partial\Omega}\|_{H^2(\Omega),L^1(\partial\Omega)}\bigr)^2 
\|u_n\|_{H^2(\Omega)}^2\Bigr)^{1/2}\\
&\leq \tilde{C}(1+\|\vec{c}\|_{\ell^\infty(L^2(\Omega))} \|u\|_{\ell^2(H^2(\Omega))}.
\end{aligned}
\]
On the other hand, for any fixed $c\in \mathcal{D}(\check{\mathbf{F}})$, as a consequence of the Lax-Milgram Lemma and elliptic regularity, the operator $\check{A}(c,u)$ is Lipschitz continuously invertible, as needed for the application of Lemma~\ref{lem:vsc_red-aao}. More precisely, for any $(\vec{f},\vec{\phi},\vec{m})\in W^*$ we have that $(\vec{f},\vec{\phi},\vec{m})=\check{A}(c,u_1)-\check{A}(c,u_2)$ iff for all $n\in\mathbb{N}$, 
$\tilde{u}_n:=u_{1,n}-u_{2,n}-m_n/|\partial\Omega|$ solves 
\[
-\Delta \tilde{u}_n +c\, \tilde{u}_n=f_n-c_n \tfrac{m_n}{|\partial\Omega|} \ \text{ in }\Omega, \quad 
\partial_\nu \tilde{u}_n=\phi_n, \ \text{ on }\partial\Omega, \quad 
\int_{\partial\Omega}\tilde{u}_{n}\, d\Gamma=0
\]
and therefore 
\[
\begin{aligned}
&\|u_1-u_2\|_{\ell^2(H^2(\Omega))}= 
\Bigl(\sum_{n\in\mathbb{N}} \|\tilde{u}_n +\tfrac{m_n}{|\partial\Omega|}\|_{H^2(\Omega)}^2\Bigr)^{1/2}\\
&\leq \Bigl(\sum_{n\in\mathbb{N}}\Bigl(
\tilde{C} \bigl(\|f_n-c_n \tfrac{m_n}{|\partial\Omega|}\|_{L^2(\Omega)}
+\|\phi_n\|_{H^{1/2}(\partial\Omega)}\bigr)
+\|\tfrac{m_n}{|\partial\Omega|}\|_{H^2(\Omega)}\Bigr)^2\Bigr)^{1/2}\\
&\leq \tilde{C} \|\vec{f}\|_{\ell^2(L^2(\Omega))}
+\tfrac{\tilde{C}}{|\partial\Omega|}\|\vec{c}\|_{\ell^\infty(L^2(\Omega))}\| \vec{m}\|_{\ell^2(\mathbb{R})}
+\tilde{C}\|\vec{\phi}\|_{\ell^2(H^{1/2}(\partial\Omega))}
+\tfrac{\sqrt{|\Omega|}}{|\partial\Omega|}\|\vec{m}\|_{\ell^2(\mathbb{R})}\\
&\leq\tilde{\tilde{C}} \|\check{A}(c,u_1)-\check{A}(c,u_2)\|_{W^*}
\end{aligned}
\]

Also Lipschitz continuity of $(I-P):X\to X$ and therefore of $P:X\to X$ is easily verified:
\begin{equation}\label{Lipschitz_P}
\|(I-P)(\vec{c}_1,u_1)-(I-P)(\vec{c}_2,u_2)\|_X
\leq\tfrac{1}{\|\vec{w}\|_{\ell^1(\mathbb{R})}} \sup_{j\in\mathbb{N}}
\sum_{k\in \mathbb{N}}w_k \|c_{1,k}-c_{2,k}\|_{\check{Q}}
\leq \|\vec{c}_1-\vec{c}_2\|_{\ell^\infty(\check{Q})}.
\end{equation}

\medskip

We make use of the following result from \cite{Diss-Weidling} that we quote in a particularized way for the convenience of the reader.
\begin{theorem}\label{thm:vsc_Schroedinger}(\cite[Theorem 6.10, 6.12 with $p=2$]{Diss-Weidling})
Let $\Omega=\mathcal{B}_1(0)\subseteq\mathbb{R}^d$, $d=3$, $\sigma_{bg}\equiv1$.
Set $\check{R}(c):=\tfrac{1}{2}\|c\|_{H^{s_1}(\Omega)}^2$ for some $s_1\in(1/2,3/2)$ and assume that $c^\dagger\in B^{s_2}_{2,\infty}$ for some $s_2\in(s_1,2s_1+1/2)$. Then there exist constants $C_1$, $C_2$, $t_0>0$ such that \eqref{vsc} holds for $\check{F}=\check{\mathbf{F}}$ with $\psi(t)=C_1(-\log(\min\{t,t_0\}+C_2))^{-2(s_2-s_1)
\frac{2s_1+1}{2s_2+1}}$ and $b=3/4$. 
\end{theorem}
This follows directly from an application of \cite[Theorem 6.12]{Diss-Weidling} with (in the notation from there) $p=q=r=2$, $s_0=-1/2$ 
to the stability result \cite[Theorem 6.10]{Diss-Weidling}.

In order to conclude a VSC for the all-at-once forward operator $\check{\mathbb{F}}$, we apply Lemma~\ref{lem:vsc_red-aao}. The purpose of using the all-at-once formulation is to allow for a bound
\begin{equation}\label{u0}
(\tfrac{1}{u_{0,j}})_{j\in\mathbb{N}}\in\ell^\infty(H^{d/2+\epsilon}(\Omega))
\end{equation}
for some $\epsilon\in(0,2-d/2]$, that will be required for a rigorous verification of the range invariance condition \eqref{rangeinvar_diff}.

Using this in the convergence results of Sections~\ref{sec:var}, \ref{sec:sm}, \ref{sec:Newton}, we obtain convergence rates for the inverse protential problem. 
\begin{corollary}
Under the conditions of Theorem~\ref{thm:vsc_Schroedinger}, with \eqref{u0}, $p=2$, $\mathcal{Q}(\hat{r}_1,\hat{r}_2)=\|\hat{r}_1-\hat{r}_2\|_X$, 
$\alpha\sim\delta^{\tilde{p}}$, $\tilde{p}\in(0,2]$, 
we obtain the convergence rate 
\[
\|\vec{c}^{\;\delta}-c^\dagger\|_{\ell^\infty(H^{s_1}(\Omega))}+\|u^\delta-S(c^\dagger)\|_{\ell^2(H^2(\Omega)}=O((-\log(\delta^2+C_2))^{-2(s_2-s_1)
\frac{2s_1+1}{(2s_2+1)p}}
\]
\begin{enumerate}
\item for $(\vec{c}^{\;\delta},u^\delta)=x_{\alpha,\beta,\eta}^\delta$ defined by \eqref{var} with \eqref{RcheckR}, \eqref{apriori};
\item for $(\vec{c}^{\;\delta},u^\delta)=x_{sm}^\delta$ defined by \eqref{sm} with \eqref{RcheckR_sm}, \eqref{gammaCbar}, \eqref{apriori};
\item for $(\vec{c}^{\;\delta},u^\delta)=x_{n_*}^\delta$ defined by \eqref{Newton} with \eqref{RcheckR_sm}, \eqref{gammaCbar}, \eqref{alpha_eta_nstar}.
\end{enumerate}
\end{corollary}
{\em Proof.}\ 
Invertibility of $r$ as defined in \eqref{r_Schroedinger}, that is, in the all-at-once setting used here, 
\begin{equation}\label{r_Schroedinger_aao}
r(\vec{c},u)_j:=
\left(\begin{array}{c}(c_j-c_{0,j})(1+\tfrac{u_j-u_{0,j}}{u_{0,j}})\\
u_j-u_{0,j}\end{array}\right), 
\quad 
r^{-1}(\hat{r}_c,\hat{r}_u)=
\left(\begin{array}{c}\hat{r}_{c,j}(1+\tfrac{\hat{r}_{u,j}}{u_{0,j}})^{-1}+c_{0,j}\\
\hat{r}_u+u_{0,j}\end{array}\right), 
\end{equation}
as well as Lipschitz continuity of $r$ and $r^{-1}$ follow from the estimate
\begin{equation}\label{Lipschitz_Schroedinger}
\begin{aligned}
&\|r(\vec{c}_1,u_1)-r(\vec{c}_2,u_2)-\bigl((\vec{c}_1,u_1)-(\vec{c}_2,u_2)\bigr)\|_X\\
&=
\|((c_{1,j}-c_{0,j})\tfrac{u_{1,j}-u_{0,j}}{u_{0,j}}
-(c_{2,j}-c_{0,j})\tfrac{u_{2,j}-u_{0,j}}{u_{0,j}})_{j\in\mathbb{N}}\|_{\ell^\infty(H^{s_1}(\Omega))}\\
&=
\sup_{j\in\mathbb{N}} 
\|\tfrac{1}{u_{0,j}}\bigl((c_{1,j}-c_{2,j})(u_{1,j}-u_{0,j})+(c_{2,j}-c_{0,j})(u_{1,j}-u_{2,j})\bigr)\|_{H^{s_1}(\Omega)}\\
&\leq C
\sup_{j\in\mathbb{N}} 
\|\tfrac{1}{u_{0,j}}\|_{H^{d/2+\epsilon}(\Omega)}
\Bigl(\|c_{1,j}-c_{2,j}\|_{H^{s_1}(\Omega)}\|u_{1,j}-u_{0,j}\|_{H^2(\Omega)}\\
&\hspace*{5cm}+\|c_{2,j}-c_{0,j}\|_{H^{s_1}(\Omega)}\|u_{1,j}-u_{2,j}\|_{H^2(\Omega)}\Bigr)
\\
&\leq \tilde{C}\|(\vec{c}_1,u_1)-(\vec{c}_2,u_2)\|_X
\bigl(\|(\vec{c}_1,u_1)-(\vec{c}_0,u_0)\|_X+\|(\vec{c}_2,u_2)-(\vec{c}_0,u_0)\|_X\bigr)\\
&\leq c_r\|(\vec{c}_1,u_1)-(\vec{c}_2,u_2)\|_X
\end{aligned}
\end{equation}
with $c_r\in(0,1)$ for $(\vec{c}_i,u_i)\in U:=\mathcal{B}_{\rho_0}^X(\vec{c}_0,u_0)$, $i\in\{1,2\}$ with $\rho_0$ small enough.
Here we have used the assumption \eqref{u0} and 
\cite[Theorem 6.1]{Johnsen1995} (see also \cite[Theorem B.22]{Diss-Weidling}) with $s_0=d/2+\epsilon$ and $p_0=p_1=p=2$ to obtain, e.g, $\|\tfrac{1}{u_{0,j}}(c_{1,j}-c_{2,j})(u_{1,j}-u_{0,j})\|_{H^{s_1}(\Omega)}\leq C 
\|c_{1,j}-c_{2,j}\|_{H^{s_1}(\Omega)}\|\tfrac{1}{u_{0,j}}(u_{1,j}-u_{0,j})\|_{H^{d/2+\epsilon}(\Omega)}$ and then the fact that $H^{d/2+\epsilon}(\Omega)$ with the pointwise multiplication is a Banach algebra as well as the embedding $\ell^2(H^2(\Omega))\subseteq\ell^\infty(H^2(\Omega))$ of the sequence spaces.

It only remains to verify the tangential cone condition for $r$, \eqref{tangcone_r}, that appears in the convergence proof of \eqref{Newton}
\[
\begin{aligned}
&\|r(\vec{c}_1,u_1)-r(\vec{c}_2,u_2)-r'(\vec{c}_2,u_2)\bigl((\vec{c}_1,u_1)-(\vec{c}_2,u_2)\bigr)\|_X\\
&=
\sup_{j\in\mathbb{N}}  
\|(c_{1,j}-c_{0,j})\tfrac{u_{1,j}-u_{0,j}}{u_{0,j}}
-(c_{2,j}-c_{0,j})\tfrac{u_{2,j}-u_{0,j}}{u_{0,j}}\\
&\hspace*{2cm}
-(c_{1,j}-c_{2,j})\tfrac{u_{2,j}-u_{0,j}}{u_{0,j}}
-(c_{2,j}-c_{0,j})\tfrac{u_{1,j}-u_{2,j}}{u_{0,j}}
\|_{H^{s_1}(\Omega)}\\
&=
\sup_{j\in\mathbb{N}} 
\|\tfrac{1}{u_{0,j}}(c_{1,j}-c_{2,j})(u_{1,j}-u_{2,j})\|_{H^{s_1}(\Omega)}\\
&\leq \tilde{C}\|(\vec{c}_1,u_1)-(\vec{c}_2,u_2)\|_X^2
\leq c_{tcr} \|r(\vec{c}_1,u_1)-r(\vec{c}_2,u_2)\|_X
\end{aligned}
\]
similarly to \eqref{Lipschitz_Schroedinger}.
\begin{flushright}
$\diamondsuit$
\end{flushright}

\medskip

Alternatively, we could apply the general relaxation from Section~\ref{sec:generalF}, by remaining with the reduced original formulation 
\[
\begin{aligned}
&\check{\mathbf{F}}:\mathcal{D}(\check{\mathbf{F}})(\subseteq \check{Q})\to Z, \quad
\check{X}=\check{Q}, \quad
Y=Z=\ell^2(L^2(\partial\Omega)), \\
&\check{\mathbf{F}}(c):= (\text{tr}^D_{\partial\Omega}u_{n})_{n\in\mathbb{N}}
\text{ where $u_{n}=S(c)_n$ and }\check{A}(c,S(c))=0
\end{aligned}
\]
with $\mathcal{D}(\check{\mathbf{F}})$ as in \eqref{DFcheck}, but defining the extended operator by 
\[
\mathbf{F}:\mathcal{D}(\check{\mathbf{F}})\times Z\to Z, \quad
X=\check{Q}\times Z, \quad
\mathbf{F}(c,z):= \mathbf{F}(c)+z
\]
and
\begin{equation}\label{r_Schroedinger_alt}
\begin{aligned}
&r(c,z)_j:=(c-c_0,\ z_j-z_{0,j}+\text{tr}^D_{\partial\Omega}v_j)\\
&\text{ where $v_j$ solves }
\left\{\begin{array}{rll}
-\Delta v_j+c\, v_j&=(c-c_0) (S(c)_j-S(c_0)_j)
 \quad &\text{ in }\Omega \\
\partial_\nu v_j&=0 &\text{ on }\partial\Omega \\ 
\int_{\partial\Omega}v_j\, d\Gamma&=0.
\end{array}\right.
\end{aligned}
\end{equation}
Since the functions $z_j$ are defined on the $d-1$ dimensional manifold $\partial\Omega$, this is a more parsimonious relaxation in the sense of augmentation of the set of unknowns. It also avoids a condition like \eqref{u0} and thus allows to stay in the reduced setting. However, this comes at the prize of the more involved definition of $r$ according to \eqref{r_Schroedinger_alt} and its inverse
\[
\begin{aligned}
&r^{-1}(\hat{r}_c,\hat{r}_z)_j:=(c_0+\hat{r}_c,\ z_{0,j}+\hat{r}_{z,j}-\text{tr}^D_{\partial\Omega}v_j)\\
&\text{ where $v_j$ solves }
\left\{\begin{array}{rll}
-\Delta v_j+(c_0+\hat{r}_c)\, v_j&= \hat{r}_c (S(c_0+\hat{r}_c)_j-S(c_0)_j)
 \quad &\text{ in }\Omega \\
\partial_\nu v_j&=0 &\text{ on }\partial\Omega \\ 
\int_{\partial\Omega}v_j\, d\Gamma&=0.
\end{array}\right.
\end{aligned}
\]
as compared to \eqref{r_Schroedinger_aao}.
\end{example}

\begin{example}[Example~\ref{ex:EIT} revisited]\label{ex:EIT_details}
The plan is similar to the one in Example~\ref{ex:Schroedinger_details} above:
Also here we make use of the function space setting detailed in \cite{Diss-Weidling} and use a variational source condition from there, but stay in the reduced setting.

We thus write the inverse problem as
\begin{equation}\label{red_EIT}
\check{\mathbf{F}}(\sigma)=y
\end{equation}
or equivalently, after extension of the parameter space and penalization according to Section~\ref{sec:generalF}, as
\begin{equation}\label{aao_EIT_extended}
\begin{aligned}
&\mathbf{F}(\sigma,z)=(0,y)^T\\
&P(\sigma,z)=0
\end{aligned}
\end{equation}
cf. \eqref{red} with 
\[
\begin{aligned}
&\check{\mathbf{F}}:\mathcal{D}(\check{\mathbf{F}})(\subseteq\check{X})\to Y, \quad 
\mathbf{F}:\mathcal{D}(\check{\mathbf{F}})\times Z\to Y, 
\\&
\mathcal{D}(\check{\mathbf{F}})= \{\sigma\in H^{s_1}(\Omega)\ : \ 0<\underline{\sigma}\leq \sigma(x)\leq \overline{\sigma}<\infty, \, x\in\Omega, \text{ and \eqref{sigmabg} holds}\}
\\
&\check{X}=H^{s_1}(\Omega), \quad
X = \check{X}\times Z, \quad
Y=Z=\ell^2(L^2(\partial\Omega))\\ 
\end{aligned}
\]
with $s_1>d/2$ (so that $\check{X}$ continuously embeds into $L^\infty(\Omega)$) and \eqref{Fred_EIT}, \eqref{P_EIT}, $\check{\mathbf{F}}(\sigma)= \mathbf{F}(\sigma,0)$.
This satisfies \eqref{rangeinvar_diff} with \eqref{r_EIT}.

Since we employ the relaxation from Section~\ref{sec:generalF}, we also have to verify the Taylor remainder estimate \eqref{Taylor} and, in case of Newton \eqref{Newton} also condition \eqref{tcr_Fcheck}.
We do so by observing that with $u_{k,j}=S(\sigma_k)_j$, $k\in\{1,2,i\}$, $i\in\{0,2\}$
\[
\begin{aligned}
&\check{\mathbf{F}}(\sigma_1)-\check{\mathbf{F}}(\sigma_2)-\check{\mathbf{F}}'(\sigma_i)=
(\text{tr}^D_{\partial\Omega}\tilde{v}_j)_{j\in\mathbb{N}} \quad \text{ where $\tilde{v}_j$ solves }
\\
&\left\{\begin{array}{rll}
-\nabla\cdot(\sigma_i\nabla \tilde{v}_j)&
=\nabla\cdot\bigl(
(\sigma_1-\sigma_i)\nabla (u_{1,j}-u_{2,j})+(\sigma_1-\sigma_2)\nabla (u_{2,j}-u_{i,j})\bigr)
 \quad &\text{ in }\Omega \\
\partial_\nu \tilde{v}_j&=0 &\text{ on }\partial\Omega \\ 
\int_{\partial\Omega}\tilde{v}_j\, d\Gamma&=0.
\end{array}\right.
\end{aligned}
\]
Thus with the Poincar\'{e}-Friedrichs type estimate
\[
\|v\|_{H^1(\Omega)} \leq C_{PF} \left(\|\nabla v\|_{L^2(\Omega)}+\left|\int_{\partial\Omega}v\, d\Gamma\right|\right), \quad
v\in H^1(\Omega)
\]
applied to the weak formulation 
\[
\begin{aligned}
&\int_\Omega \sigma_i\nabla \tilde{v}_j\cdot \nabla w\, dx =
-\int_\Omega \bigl((\sigma_1-\sigma_i)\nabla (u_{1,j}-u_{2,j})+(\sigma_1-\sigma_2)\nabla (u_{2,j}-u_{i,j})\bigr)\cdot \nabla w\, dx, \\  
&\text{for all }w\in H^1(\Omega) \text{ with }\int_{\partial\Omega}w\, d\Gamma =0,
\end{aligned}
\]
setting $w=\tilde{v}_j$ and applying the Cauchy-Schwarz inequality, we obtain
\[
\begin{aligned}
&\|\check{\mathbf{F}}(\sigma_1)-\check{\mathbf{F}}(\sigma_2)-\check{\mathbf{F}}'(\sigma_i)\|_{L^2(\partial \Omega)}\\
&\leq \tfrac{1}{\underline{\sigma}} C_{PF} \|\text{tr}^D_{\partial\Omega}\|_{H^1(\Omega),L^2(\partial\Omega)}
\|(\sigma_1-\sigma_i)\nabla (u_{1,j}-u_{2,j})+(\sigma_1-\sigma_2)\nabla (u_{2,j}-u_{i,j})\|_{L^2(\Omega)}
\\
&\leq \tfrac{1}{\underline{\sigma}} C_{PF} \|\text{tr}^D_{\partial\Omega}\|_{H^1(\Omega),L^2(\partial\Omega)}
\Bigl(\|\sigma_1-\sigma_i\|_{L^\infty(\Omega)} \|u_{1,j}-u_{2,j}\|_{H^1(\Omega)}\\
&\hspace*{4.5cm}+\|\sigma_1-\sigma_2\|_{L^\infty(\Omega)} \|u_{2,j}-u_{i,j}\|_{H^1(\Omega)}\Bigr)\\
&\leq \tilde{C} \Bigl(\|\sigma_1-\sigma_i\|_{\check{X}}+\|\sigma_2-\sigma_i\|_{\check{X}}\Bigr)
\|\sigma_1-\sigma_2\|,
\end{aligned}
\]
due to Lipschitz continuity of $S:L^\infty(\Omega)\to H^1(\Omega)$. \\
Thus \eqref{Taylor} and \eqref{tcr_Fcheck} hold on a sufficiently small neighborhood $U$ of $\sigma^\dagger$.

\begin{theorem}\label{thm:vsc_EIT}(\cite[Theorem 6.15 with $p=2$]{Diss-Weidling})
Let $\Omega=\mathcal{B}_1(0)\subseteq\mathbb{R}^d$, $d=3$, $\sigma_{bg}\equiv1$. 
Set $\check{R}(\sigma):=\tfrac{1}{2}\|\sigma\|_{H^{s_1}(\Omega)}^2$ for some $s_1\in(5/2,7/2)$ and assume that $\sigma^\dagger\in B^{s_2}_{2,\infty}$ for some $s_2\in(s_1,2s_1-1)$. Then there exist constants $C$, $t_0>0$
such that \eqref{vsc} holds with $\psi(t)=C(-\log(\min\{t,t_0\}))^{-(s_2-s_1)\frac{2s_1-3}{s_1-1}}$ and $b=3/4$. 
\end{theorem}

Combining this with the convergence results of Sections~\ref{sec:var}, \ref{sec:sm}, \ref{sec:Newton}, we arrive at the following convergence rates result for EIT. 
\begin{corollary}
Under the conditions of Theorem~\ref{thm:vsc_EIT}, with $p=2$, $\mathcal{Q}(\hat{r}_1,\hat{r}_2)=\|\hat{r}_1-\hat{r}_2\|_X$, 
$\alpha\sim\delta^{\tilde{p}}$, $\tilde{p}\in(0,2]$, 
we obtain the convergence rate 
\[
\|\sigma^\delta-\sigma^\dagger\|_{H^{s_1}(\Omega)}+\|z^\delta\|_{\ell^2(L^2(\partial\Omega))}=O((-\log(\delta^p+C_2))^{-(s_2-s_1)\frac{2s_1-3}{(s_1-1)p}}
\]
\begin{enumerate}
\item for $(\sigma^\delta,z^\delta)=x_{\alpha,\beta,\eta}^\delta$ defined by \eqref{var_eta} with \eqref{RcheckR}, \eqref{apriori};
\item for $(\sigma^\delta,z^\delta)=x_{sm}^\delta$ defined by \eqref{sm_eta} with \eqref{RcheckR_sm}, \eqref{gammaCbar}, \eqref{apriori};
\item for $(\sigma^\delta,z^\delta)=x_{n_*}^\delta$ defined by \eqref{Newton_eta} with \eqref{RcheckR_sm}, \eqref{gammaCbar}, \eqref{alpha_eta_nstar}.
\end{enumerate}
\end{corollary}
\end{example}

\renewcommand\appendixname{Appendix}
\begin{appendices}
\section{Appendix}
\subsection{Proof of Lemma~\ref{minLemma}}
{\em Proof.}\ 
Splitting the $\res $ term into a convex combination $\res=\lambda \res+(1-\lambda) \res$ we obtain from \eqref{minLemma}
\[
\lambda c_1\res  + c_2\alpha\, \err 
\leq d + \alpha\Bigl(-\frac{(1-\lambda) c_1}{C_3\alpha}\, C_3\res-(-\psi)(C_3\res)\Bigr)
\leq d + \alpha (-\psi)^*\Bigl(-\frac{(1-\lambda) c_1}{C_3\alpha}\Bigr)
\]
Setting $\lambda=0$ yields the estimate on $\err$ and setting $\lambda=\frac12$ the one on $\res$.
\begin{flushright}
$\diamondsuit$
\end{flushright}

\subsection{Convergence of \eqref{var} without VSC}
under the following conditions.
\begin{assumption}\label{ass:conv}
$r$ is injective and there exist tolopogies $\widetilde{\mathcal{T}}$, $\widetilde{\widetilde{\mathcal{T}}}$, on $\widetilde{X}$ and $\mathcal{T}$ on $X$ such that
\begin{enumerate}
\item[(i)]
sublevel sets of $\mathcal{R}$ are $\widetilde{\mathcal{T}}$ compact; 
\item[(ii)] For any two sequences $(\hat{r}^1_k)_{k\in\mathbb{N}}$, $(\hat{r}^2_k)_{k\in\mathbb{N}}$
\[
\mathcal{Q}(\hat{r}^1_k,\hat{r}^2_k)\to0\text{ and }
\hat{r}^2_k\stackrel{\widetilde{\mathcal{T}}}{\longrightarrow}\hat{r}^2
\ \Rightarrow \ \text{ there exists a subsequence }
\hat{r}^1_{k_\ell}\stackrel{\widetilde{\widetilde{\mathcal{T}}}}{\longrightarrow}\hat{r}^2;
\]
\item[(iii)] 
$r^{-1}:r(U)\to U$ is $\widetilde{\widetilde{\mathcal{T}}}-\mathcal{T}$ continuous;
\item[(iv)] $\|K\cdot\|_Y$ is $\widetilde{\mathcal{T}}$ lower semicontinuous;
\item[(v)] $\mathcal{P}$ is $\mathcal{T}$ lower semicontinuous.
\end{enumerate}
\end{assumption}
Moreover, the regularization and tolerance parameters $\alpha=\alpha(\delta)$, $\beta=\beta(\delta)$, $\eta(\delta)$ are assumed to satisfy 
\begin{equation}\label{alphabetaeta}
\alpha(\delta)\to0, \qquad
\frac{\delta^p+\eta(\delta)}{\alpha(\delta)}\to0, \qquad
\frac{\delta^p+\eta(\delta)}{\beta(\delta)}\to0, \qquad\text{ as }\delta\to0.
\end{equation}
The a priori choice \eqref{apriori} is compatible with \eqref{alphabetaeta}, and so is the particular one in the logarithmic case indicated in Remark~\ref{rem:alphachoice} when slightly restricted to $\alpha(\delta)\sim\delta^{\tilde{p}}$ with $\tilde{p}\in(0,p)$.

\begin{proposition}\label{prop:conv}
Under condition \eqref{rangeinvar_diff} and Assumption~\ref{ass:conv}, the approximations $x_{\alpha,\beta,\eta}^\delta$ defined by \eqref{var_eta} with \eqref{alphabetaeta} converge $\mathcal{T}$ subseqentially to a solution of \eqref{FP}, 
that is, every subsequence of  $(\hat{r}_{\alpha,\beta,\eta}^\delta,x_{\alpha,\beta,\eta}^\delta)_{\delta>0}$ has a $\widetilde{\mathcal{T}}\times\mathcal{T}$ convergent subsequence and the limit $(\hat{r},x)$ of every $\widetilde{\mathcal{T}}\times\mathcal{T}$ convergent subsequence solves \eqref{FP}. 
If the $\mathcal{R}$ minimizing solution $(\hat{r}^\dagger,x^\dagger)$ of 
\eqref{FP}\footnote{i.e., $(\hat{r}^\dagger,x^\dagger)$ solving \eqref{FP} such that 
$\forall (\hat{r},x)\mbox{ solving \eqref{FP}}\,: \  \mathcal{R}(\hat{r}^\dagger,x^\dagger)\leq \mathcal{R}(\hat{r},x)$}
is unique and 
then $\hat{r}^\delta\stackrel{\widetilde{\mathcal{T}}}{\longrightarrow}r(x^\dagger)$ and $x^\delta\stackrel{\mathcal{T}}{\longrightarrow}x^\dagger$.
\end{proposition}
{\em Proof.}\ 
For an arbitrary sequence $\delta_n\to0$, abbreviating $(\hat{r}_n,x_n):=(\hat{r}_{\alpha_n\beta_n\eta_n}^{\delta_n},x_{\alpha_n\beta_n\eta_n}^{\delta_n}$, from \eqref{min0_var} (with $(2+C_\eta)\delta^p$ replaced by the more general bound $2\delta^p+\eta(\delta)$ and \eqref{alphabetaeta} we obtain
\[
\begin{aligned}
(a)&\limsup_{n\to\infty} \mathcal{R}(\hat{r}_n) \leq \mathcal{R}(r(x^\dagger))\\
(b)&\limsup_{n\to\infty} \mathcal{Q}(r(x_n),\hat{r}_n)\ =0\\
(c)&\limsup_{n\to\infty} \|K(\hat{r}_n-r(x^\dagger))\|_Y\ =0\\
(d)&\limsup_{n\to\infty} \|P(x_n)\| =0.
\end{aligned}
\]
From (i) and (a) we conclude existence of a subsequence $\hat{r}_{n_k}$ an and element $\hat{r}\in\widetilde{X}$ with $\mathcal{R}(\hat{r})\leq\mathcal{R}(r(x^\dagger))$ such that $\hat{r}_{n_k}\stackrel{\widetilde{\mathcal{T}}}\to\hat{r}$. This together with (b), (ii) yields existence of a subsequence such that $r(x_{n_{k_\ell}})\stackrel{\tilde{\tilde{\mathcal{T}}}}\to \hat{r}$.
Now by (iii) we obtain $x_{n_{k_\ell}}\stackrel{\mathcal{T}}\to x$ and $r(x)=\hat{r}$.
Finally, (c), (iv) imply $K(r(x))=K(r(x^\dagger))$ and (d), (v) imply $P(x)=0$.
\begin{flushright}
$\diamondsuit$
\end{flushright}

\subsection{Verification of \eqref{Cpsi} for H\"older or logarithmic source conditions \eqref{psiHoelderlog}}\label{sec:ratescases}
\begin{lemma}\label{lem:Cpsi}
The functions defined in \eqref{psiHoelderlog} satisfy
\[
\psi(a+b)\leq \psi(a)+\psi(b), \text{ for all }a,b>0
\]
\end{lemma}
{\em Proof.}\ 
In the H\"older case, with $m:=1/\mu$, $x=(b/a)^\mu$ the claim is equivalent to 
\[
\phi(x):=\frac{1+x^m}{(1+x)^m}\leq 1, \text{ for all }x>0
\]
which in case of an integer $m$ immediately follows from the binomial theorem, otherwise from the fact that
$\lim_{x\searrow0} \phi(x)=1$, $\lim_{x\nearrow\infty} \phi(x)=1$, 
$\phi'(x)=(1+x)^{-2m}(mx^{m-1}(1+x)^m-(1+x^m)m(1+x)^{m-1}=
m (1+x)^{-m-1}(x^{m-1}-1)=0$ iff $x=1$ and $\phi(1)=2^{1-m}\leq 1$, since $m=1/\mu\geq1$.

In the logarithmic case for simplicity of exposition we just consider $\nu=1$, $t_0=1/e$. 
With $\alpha:=(-\log a)^{-1}$, $\beta:=(-\log b)^{-1}$, the claim is equivalent to 
\[
\phi(\beta):=\frac{-\log\left(e^{-1/\alpha}+e^{-1/\beta}\right)}{\alpha+\beta}\leq 1, \text{ for all }\beta>0
\]
and any fixed $\alpha>0$.
It is readily checked that $\lim_{\beta\searrow0} \phi(\beta)=1$, $\lim_{\beta\nearrow\infty} \phi(\beta)=0$ and  
$\phi'(\beta_*)=(\alpha+\beta_*)^{-2}\Bigl(-\left(e^{-1/\alpha}+e^{-1/\beta_*}\right)^{-1}e^{-1/\beta_*}\tfrac{\alpha+\beta_*}{\beta_*^2}+\log\left(e^{-1/\alpha}+e^{-1/\beta_*}\right)\Bigr)=0$ iff
$\log\left(e^{-1/\alpha}+e^{-1/\beta_*}\right)=
\left(e^{-1/\alpha}+e^{-1/\beta_*}\right)^{-1}e^{-1/\beta_*}\tfrac{\alpha+\beta_*}{\beta_*^2}\geq0$
which inserted into $\phi$ yields
$\phi(\beta_*)=-\left(e^{-1/\alpha}+e^{-1/\beta_*}\right)^{-1}e^{-1/\beta_*} \, \beta_*^{-2}
=-\frac{b_*}{a+b_*}(-\log b_*)^{-2}\leq0$.
\begin{flushright}
$\diamondsuit$
\end{flushright}
\end{appendices}

\end{document}